\documentclass{lmcs} 
\pdfoutput=1

\usepackage{lastpage}
\lmcsdoi{17}{2}{7}
\lmcsheading{}{\pageref{LastPage}}{}{}%
{May~01,~2019}{Apr.~20,~2021}{}

\usepackage[utf8]{inputenc}
\usepackage{booktabs}
\usepackage{dcolumn}
\newcolumntype{d}[1]{D{.}{.}{#1}}


\theoremstyle{plain} 

\RequirePackage{bussproofs}

\newcommand{\C}[1]{\mathcal{#1}}

\newcommand{\bigland}{\mathop{\hbox{$\bigwedge$\kern-.6em$\bigwedge$}}}
\newcommand{\biglor}{\mathop{\hbox{$\bigvee$\kern-.6em$\bigvee$}}}
\newcommand{\ex}{\exists} 

\newcommand{\exE}{{\ex}^{-}}

\newcommand{\constr}{\mathrm{C}}
\newcommand{\rec}{\C{R}}
\newcommand{\constrD}{\mathrm{D}} 
\newcommand{\constrFin}{\mathrm{Fin}} 

\newcommand{\constrIntN}{\mathrm{IntN}}
\newcommand{\constrIntP}{\mathrm{IntP}}

\newcommand{\constrLr}{\mathrm{Lr}} 

\newcommand{\constrU}{\mathrm{U}} 
\newcommand{\false}{\mathsf{f\kern-0.14em f}}
\newcommand{\one}{\mathrm{1}}
\newcommand{\suc}{S}
\newcommand{\termSumIntroLeft}{\mathrm{InL}}
\newcommand{\termSumIntroRight}{\mathrm{InR}}
\newcommand{\termUnit}{\mathrm{Dummy}}

\newcommand{\true}{\mathsf{t\kern-0.14em t}}

\newcommand{\typeB}{\mathbb{B}}
\newcommand{\typeBin}{\mathbb{Y}}
\newcommand{\typeIntv}{\mathbb{S}} 
\newcommand{\typeL}[1]{\mathbb{L}(#1)} 
\newcommand{\typeList}{\mathbb{L}} 
\newcommand{\typeN}{\mathbb{N}}
\newcommand{\typeG}{\mathbb{G}}
\newcommand{\typeH}{\mathbb{H}}
\newcommand{\typeI}[1]{\mathbb{I}(#1)} 
\newcommand{\typeP}{\mathbb{P}}
\newcommand{\typeProd}{\times}
\newcommand{\typeQ}{\mathbb{Q}}
\newcommand{\typeR}{\mathbb{R}}
\newcommand{\typeS}[1]{\mathbb{S}(#1)} 
\newcommand{\typeSd}{\mathbb{D}}
\newcommand{\typeSds}{\mathbb{S}} 

\newcommand{\typeSum}{+}
\newcommand{\typeTo}{\to}
\newcommand{\typeUnit}{\mathbb{U}}

\newcommand{\typeZ}{\mathbb{Z}}

\newcommand{\CoG}{{}^{\mathrm{co}}\!G}
\newcommand{\CoH}{{}^{\mathrm{co}}\!H}

\newcommand{\CoI}{{}^{\mathrm{co}}\!I}
\newcommand{\CoTotal}{{}^{\mathrm{co}}\Total}

\newcommand{\TCF}{\mathrm{TCF}} 

\newcommand{\extrTer}[1]{\mathrm{et}(#1)}
\newcommand{\extrTy}[1]{\tau(#1)}

\newcommand{\BNFdef}{\mathtt{\; ::= \;}}
\newcommand{\BNFor}{\mid}

\newcommand{\consistent}{\uparrow}
\newcommand{\corec}{{}^{\mathrm{co}}\C{R}}

\newcommand{\defeq}{:=}

\newcommand{\destr}{\C{D}}
\newcommand{\entails}{\vdash}
\newcommand{\eqd}[2]{#1 \equiv #2}

\newcommand{\inquotes}[1]{``#1''}
\newcommand{\lev}[1]{\mathrm{lev}(#1)}
\newcommand{\mr}{\mathrel{\mathbf{r}}}
\newcommand{\mrind}{\textbf{r}}

\newcommand{\ncind}[1]{#1^{\mathrm{nc}}}
\newcommand{\nil}{[]}
\newcommand{\pair}[2]{\langle #1 , #2 \rangle}
\newcommand{\set}[2]{\{\,#1\mid#2\,\}}
\newcommand{\shG}{\mathrm{shG}}
\newcommand{\shH}{\mathrm{shH}}
\newcommand{\Total}{T}

\theoremstyle{definition}
\newtheorem*{definition*}{Definition}
\newtheorem*{examples*}{Examples}

\begin{document}

\title{Logic for exact real arithmetic}

\author[H.~Schwichtenberg]{Helmut Schwichtenberg\rsuper{a}}	
\address{\lsuper{a}Mathematisches Institut, LMU, M\"unchen}	
\email{schwicht@math.lmu.de}  

\author[F.Wiesnet]{Franziskus Wiesnet\rsuper{b}}	
\address{\lsuper{b}Mathematisches Institut, LMU, M\"unchen and Dipartimento di Matematica, Università degli Studi di Trento}	
\email{franziskus.wiesnet@unitn.it}  
\thanks{This project has received funding from the European Union's
  2020 research and innovation programme under the Marie
  Sk{\l}odowska-Curie grant agreement No 731143.  The second author is a
  Marie Sk{\l}odowska-Curie fellow of the Istituto Nazionale di Alta
  Matematica}





\begin{abstract}
  \noindent   Continuing earlier work of the first author with U.~Berger,
  K.~Miyamoto and H.~Tsuiki, it is shown how a division algorithm for
  real numbers given as a stream of signed digits can be extracted
  from an appropriate formal proof.  The property of being a real
  number represented as a stream is formulated by means of
  coinductively defined predicates, and formal proofs involve
  coinduction.  The proof assistant Minlog is used to generate the
  formal proofs and extract their computational content as terms of
  the underlying theory, a form of type theory for finite or infinite
  data.  Some experiments with running the extracted term are
  described, after its translation to Haskell.
\end{abstract}

\maketitle
Real numbers in the exact (as opposed to floating-point) sense can be
defined as Cauchy sequences (of rationals, with modulus).  However,
for computational purposes it is better to see them as coded by
\inquotes{streams} of signed digits $\{1,0,-1\}$.  A variant stream
representation is the so-called \inquotes{binary reflected} or
Gray-code \cite{Gianantonio99, Tsuiki02} explained below.  Apart from
being practically more useful, the stream view turns real numbers into
\inquotes{infinite data} and hence objects of type level $0$.  As a
consequence the type level of other concepts in constructive analysis
\cite{Bishop67} is lowered by one, which simplifies matters
considerably.

Our overall goal is to obtain formally verified algorithms (given by
terms in our language) operating on stream represented real numbers.
Given an informal idea of how the algorithm should work, there are two
methods how this can be achieved.
\begin{enumerate}
\item[(I)] Formulate (using corecursion) the algorithm in the term
  language of a suitable theory, and then formally prove that this
  term satifies the specification;
\item[(II)] Find a formal existence proof $M$ (using coinduction) for
  the object the algorithm is supposed to return.  Then apply a proof
  theoretic method (\inquotes{realizability}) to extract $M$'s
  computational content as a term (involving corecursion) in the term
  language of the underlying theory.  The verification is done by a
  formal soundness proof of the realizability interpretation.  The
  extraction of the computational content and the verification are
  automatic.
\end{enumerate}
A general advantage of (II) over (I) is that one does not need to
begin with a detailed formulation of the algorithm, but instead can
stay on a more abstract level when proving the (existential)
specification.  In mathematics we know how to organize proofs, for
instance by splitting them into lemmas or on occasion make use of more
abstract concepts.  In short, mathematical experience can help to find
a well-structured algorithmic solution.

Method (I) was employed in \cite{CiaffaglioneGianantonio06} using Coq,
and method (II) in \cite{Berger09, MiyamotoSchwichtenberg15,
  BergerMiyamotoSchwichtenbergTsuiki16} using Minlog.
The present paper continues previous work \cite{MiyamotoSchwichtenberg15,
  BergerMiyamotoSchwichtenbergTsuiki16} by a case study on division.
It is shown how a division algorithm for real numbers represented as
either streams of signed digits or else in binary reflected form (Gray
code) can be extracted from an appropriate formal proof dealing with
concrete real numbers (Cauchy sequences of rationals, with moduli); in
\cite{MiyamotoSchwichtenberg15, BergerMiyamotoSchwichtenbergTsuiki16}
real numbers were treated axiomatically.  The reason for switching to
concrete real numbers is the desire to have reliable programs, which
are less studied for real number computation.  The method described in
this paper achieves reliability by providing an (automatically
generated) formal proof that the extracted program meets its
specification.  These proofs are in a formal theory whose consistency
is guaranteed by a concrete model \cite{Scott70, Ershov77}.  There is
no need to rely on an axiom system for real numbers as in
\cite{MiyamotoSchwichtenberg15, BergerMiyamotoSchwichtenbergTsuiki16}.

We will work with constructive existence proofs in the style of
\cite{Bishop67}, but in such a way that we can switch on and off the
availability of input data for the constructions implicit in the proof
\cite{Berger93a}.  In the present context this will be applied to real
numbers as input data: we do not want to make use of the Cauchy
sequence for the constructions to be done, but only the computational
content of an appropriate coinductive predicate to which the real
number is supposed to belong.  We consider division of real numbers
as a non-trivial case study; it has been dealt with in
\cite{CiaffaglioneGianantonio06} using method (I).  Based on the
algorithmic idea in \cite{CiaffaglioneGianantonio06}, we employ method
(II) to extract signed digit and Gray code based stream algorithms for
division from proofs that the reals are closed under division (under
some obvious restrictions).  As a benefit from the necessary
organization into a sequence of lemmas we obtain a relatively easy
analysis of the \inquotes{look-ahead}, i.e., how far we have to look
into the argument streams to obtain the $n$-th digit of the result
stream.

The paper is organized as follows.  Section~\ref{S:Prelim} collects
some background material.  Section~\ref{SS:Model} describes (mainly by
examples) the base type objects of the Scott-Ershov \cite{Scott70,
  Ershov77} model of partial continuous functionals, which is the
intended model of the theory $\TCF$ we are going to use.  It is a form
of type theory suitable for reasoning with finite and infinite data.
This theory is described in Section~\ref{SS:TCF}.  The notion of
computational content of proofs is introduced in
Section~\ref{SS:CompContent}.  Section~\ref{SS:Reals} recalls the
constructive theory of real numbers of Bishop \cite[Ch.2]{Bishop67}.
In Section~\ref{S:Div} division for stream-represented real numbers is
studied, both for signed and binary reflected digit streams.
Formalizations of the relevant proofs are described, including the
terms (programs) extracted from them.  Numerical experiments after
translating these terms into Scheme or Haskell programs are discussed.
Section~\ref{S:Sound} recalls the notion of realizability in
Section~\ref{SS:Realizers}, sketches a proof of the general soundness
theorem in Section~\ref{SS:Sound}, and its formalization for real
division algorithms in Section~\ref{SS:FSound}.
Section~\ref{S:Conclusion} concludes.

\section{Preliminaries}
\label{S:Prelim}
We informally describe a formal theory $\TCF$ \cite[Section
  7.1]{SchwichtenbergWainer12} which will be used to carry out proofs
on real numbers and also functions and predicates on real numbers.  In
contrast to \cite{MiyamotoSchwichtenberg15,
  BergerMiyamotoSchwichtenbergTsuiki16} where an axiomatic approach
was used we now work with concrete real numbers: Cauchy sequences of
rationals, with moduli.  The reason for this switch is that we aim at
full reliability of the soundness proofs we will (automatically)
generate (see Section~\ref{S:Sound}).  The axioms the formal theory is
based upon, and hence all propositions proved, are valid in a model
and therefore correct.

The base type data of the theory are free algebras given by their
constructors, for instance lists of signed digits.  Functions
(\inquotes{program constants}) are defined by defining equations
(\inquotes{computation rules}) whose left hand sides are
non-overlapping constructor patterns.  Termination is \emph{not}
required, which makes it possible the have the corecursion operator
(see below) as an ordinary constant in the language.  Predicates are
least and greatest fixed points of clauses..  Examples are totality
and cototality of lists of signed digits.  The only axioms are the
defining equations of the functions and the introduction and
elimination axioms for the (co)inductively defined constants, which
state their fixed point properties.

An important predicate is equality.  Clearly we have the inductively
defined Leibniz equality and the decidable equality for finitary
algebras like the natural numbers.  Equality on the (concrete) reals
is more delicate.  It is defined from an inductive $\le$-predicate
which refers to the underlying Cauchy sequences and their moduli.  For
all functions and predicates on the reals that we consider it is
necessary to prove compatibility with real equality.  The reason for
this extension of the axiomatic treatment is to obtain formal proofs
which rest only on the fixed point axioms and in this sense are
correct: the propositions proved are valid in the concrete model of
$\TCF$.

\subsection{Objects}
\label{SS:Model}
Base type objects of our theory are determined by the constructors of
(free) algebras.  More precisely, let $\alpha, \beta, \xi$ be type
variables.

\begin{definition*}
  \emph{Constructor types}\index{constructor type} $\kappa$ have the
  form
  \begin{equation*}
    \vec{\alpha} \typeTo (\xi)_{i<n} \typeTo \xi
  \end{equation*}
  with all type variables $\alpha_i$ distinct from each other and from
  $\xi$.  Iterated arrows are understood as associated to the right.
  An argument type of a constructor type is called a \emph{parameter}%
  \index{parameter!argument type} argument type if it is different from
  $\xi$, and a \emph{recursive} argument type%
  \index{recursive!argument type} otherwise.  A constructor type
  $\kappa$ is \emph{nullary}%
  \index{nullary!argument type} if it has no recursive argument types.
  We call
  \begin{equation*}
    \iota \defeq \mu_\xi \vec{\kappa}
    \quad \hbox{($\vec{\kappa}$ not empty)}
  \end{equation*}
  an \emph{algebra form}%
  \index{algebra form} (\inquotes{form} since type parameters may occur).  
  An algebra form is \emph{non-recursive}%
  \index{algebra form!non-recursive}\index{algebra form!explicit} (or
  \emph{explicit}) if it does not have recursive argument types.
\end{definition*}

\begin{examples*}
  We list some algebra forms without parameters, with standard names
  for the constructors added to each constructor type.
  \begin{alignat*}{3}
    &\typeUnit &&\defeq \mu_\xi( \termUnit \colon \xi)
    &&\qquad \hbox{(unit\index{unit})},
    \\
    &\typeB &&\defeq \mu_\xi (\true \colon \xi, \false \colon \xi)
    &&\qquad \hbox{(booleans\index{boolean})},
    \\
    &\typeSd
    &&\defeq \mu_\xi
    (\mathrm{SdR} \colon \xi, \mathrm{SdM} \colon \xi, \mathrm{SdL} \colon \xi)
    &&\qquad
    \hbox{(signed digits, for $1$, $0$, $-1$)},
    \\
    &\typeN &&\defeq \mu_\xi (0 \colon \xi, \suc \colon \xi \typeTo \xi) 
    &&\qquad
    \hbox{(natural numbers%
      \index{number}, unary%
      \index{number!unary})},
    \\
    &\typeP &&\defeq \mu_\xi (\one \colon \xi, \suc_0 \colon \xi \typeTo \xi,
    \suc_1 \colon \xi \typeTo \xi)
    &&\qquad \hbox{(positive numbers%
      \index{number!positive}, binary%
      \index{number!binary})},
    \\
    &\typeZ &&\defeq \mu_\xi (
    0_\typeZ \colon \xi,
    \constrIntP \colon \typeP \typeTo \xi,
    \constrIntN \colon \typeP \typeTo \xi)
    &&\qquad \hbox{(integers)},
    \\
    &\typeBin &&\defeq \mu_\xi (- \colon \xi,
    \mathrm{Branch} \colon \xi \typeTo \xi \typeTo \xi)
    &&\qquad \hbox{(binary trees%
      \index{binary tree}\index{tree!binary})}.
  \end{alignat*}
  Algebra forms with type parameters are
  \begin{alignat*}{3}
    &\typeI{\alpha} &&\defeq \mu_\xi (\mathrm{Id} \colon \alpha \typeTo \xi)
    &&\qquad \hbox{(identity\index{identity})},
    \\
    &\typeL{\alpha} &&\defeq \mu_\xi (\mathrm{Nil} \colon \xi,
    \mathrm{Cons} \colon \alpha \typeTo \xi \typeTo \xi)
    &&\qquad \hbox{(lists\index{list})},
    \\
    &\typeS{\alpha} &&\defeq \mu_\xi (\mathrm{SCons} \colon
    \alpha \typeTo \xi \typeTo \xi)
    &&\qquad \hbox{(streams\index{stream})},
    \\
    &\alpha \typeProd \beta &&\defeq \mu_\xi (
    \mathrm{Pair} \colon \alpha \typeTo \beta \typeTo \xi)
    &&\qquad \hbox{(product\index{product})},
    \\
    &\alpha \typeSum \beta &&\defeq \mu_\xi (
    \termSumIntroLeft \colon \alpha \typeTo \xi,
    \termSumIntroRight \colon \beta \typeTo \xi)
    &&\qquad \hbox{(sum\index{sum})}.
  \end{alignat*}
  The default name for the $i$-th constructor of an algebra form is
  $\constr_i$.
\end{examples*}

\begin{definition*}[Type\index{type}]
  \begin{equation*}
    \rho, \sigma, \tau \BNFdef
    \alpha \BNFor \iota(\vec{\rho}\,) \BNFor \tau \typeTo \sigma,
  \end{equation*}
  where $\iota$ is an algebra form with $\vec{\alpha}$ its parameter type
  variables, and $\iota(\vec{\rho}\,)$ the result of substituting the
  (already generated) types $\vec{\rho}$ for $\vec{\alpha}$.  Types of
  the form $\iota(\vec{\rho}\,)$ are called
  \emph{algebras}\index{algebra}.  An algebra is \emph{closed}%
  \index{algebra!closed} if it has no type variables.  
  The \emph{level}%
  \index{type!level of}\index{level!of a type} of a type is defined by
  \begin{align*}
    \lev{\alpha} &\defeq 0,
    \\
    \lev{\iota(\vec{\rho}\,)} &\defeq \max(\lev{\vec{\rho}\,}),
    \\
    \lev{\tau \typeTo \sigma} &\defeq \max( \lev{\sigma}, 1+\lev{\tau}).
  \end{align*}
  \emph{Base} types%
  \index{type!base} are types of level $0$, and a \emph{higher}
  type\index{type!higher} has level at least $1$.  
\end{definition*}

\begin{examples*}
  1.  $\typeL{\alpha}$, $\typeL{\typeL{\alpha}}$, $\alpha
  \typeProd \beta$ are algebras.

  2.  $\typeL{\typeL{\typeN}}$, $\typeN \typeSum \typeB$, $\typeQ \defeq
  \typeZ \typeProd \typeP$ are closed base types.

  3.  $(\typeN \typeTo \typeQ) \typeProd (\typeP \typeTo
  \typeN)$ is a closed algebra of level $1$.
\end{examples*}

There can be many equivalent ways to define a particular type.  For
instance, we could take $\typeUnit \typeSum \typeUnit$ to be the type
of booleans, $\typeL{\typeUnit}$ to be the type of natural numbers,
and $\typeL{\typeB}$ to be the type of positive binary numbers.

We introduce the base type objects by example only, for $\typeN$,
$\typeBin$, $\typeL{\typeB}$; cf.\ \cite[Section
  6.1.4]{SchwichtenbergWainer12} for the general case.  These objects
are built from \emph{tokens}, which are type correct constructor
expressions possibly containing the special symbol $\ast$ which
carries no information.  Examples in the algebra $\typeBin$ of binary
trees are
\begin{equation*}
  \mathrm{Branch}(-, \mathrm{Branch}(-,-))
  \qquad
  \qquad
  \mathrm{Branch}(-, \mathrm{Branch}(\ast,-))
\end{equation*}
or pictorially
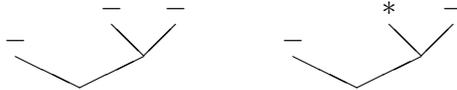
\begin{equation*}
  \begin{picture}(60,36)
    \put(0,18){\makebox(0,0){$-$}}
    \put(36,30){\makebox(0,0){$-$}}
    \put(60,30){\makebox(0,0){$-$}}
    
    \put(0,12){\line(2,-1){24}}
    \put(24,0){\line(2,1){24}}
    \put(48,12){\line(-1,1){12}}
    \put(48,12){\line(1,1){12}}
  \end{picture}
  \qquad
  \qquad
  \begin{picture}(60,36)
    \put(0,18){\makebox(0,0){$-$}}
    \put(36,30){\makebox(0,0){$\ast$}}
    \put(60,30){\makebox(0,0){$-$}}
    
    \put(0,12){\line(2,-1){24}}
    \put(24,0){\line(2,1){24}}
    \put(48,12){\line(-1,1){12}}
    \put(48,12){\line(1,1){12}}
  \end{picture}
\end{equation*}
For tokens we have obvious concepts called consistency ($\consistent$)
and entailment ($\entails$).  For instance
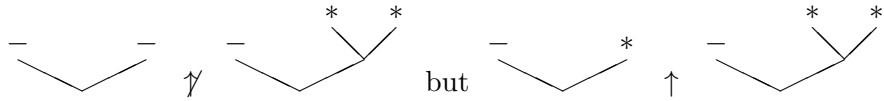
\begin{equation*}
  \begin{picture}(48,36)
    \put(0,18){\makebox(0,0){$-$}}
    \put(48,18){\makebox(0,0){$-$}}
     
    \put(0,12){\line(2,-1){24}}
    \put(24,0){\line(2,1){24}}
  \end{picture}
  \quad
  \not\consistent
  \quad
  \begin{picture}(60,36)
    \put(0,18){\makebox(0,0){$-$}}
    \put(36,30){\makebox(0,0){$\ast$}}
    \put(60,30){\makebox(0,0){$\ast$}}
    
    \put(0,12){\line(2,-1){24}}
    \put(24,0){\line(2,1){24}}
    \put(48,12){\line(-1,1){12}}
    \put(48,12){\line(1,1){12}}
  \end{picture}
  \quad\hbox{but}\quad
  \begin{picture}(48,36)
    \put(0,18){\makebox(0,0){$-$}}
    \put(48,18){\makebox(0,0){$\ast$}}
     
    \put(0,12){\line(2,-1){24}}
    \put(24,0){\line(2,1){24}}
  \end{picture}
  \quad
  \consistent
  \quad
  \begin{picture}(60,36)
    \put(0,18){\makebox(0,0){$-$}}
    \put(36,30){\makebox(0,0){$\ast$}}
    \put(60,30){\makebox(0,0){$\ast$}}
    
    \put(0,12){\line(2,-1){24}}
    \put(24,0){\line(2,1){24}}
    \put(48,12){\line(-1,1){12}}
    \put(48,12){\line(1,1){12}}
  \end{picture}
\end{equation*}
Examples for entailment in $\typeBin$ are
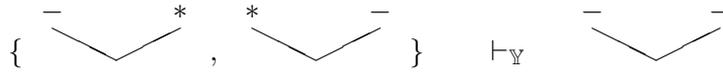
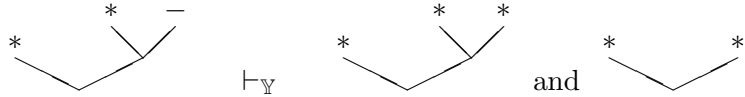
\begin{equation*}
  \{\quad
  \begin{picture}(48,36)
    \put(0,18){\makebox(0,0){$-$}}
    \put(48,18){\makebox(0,0){$\ast$}}
     
    \put(0,12){\line(2,-1){24}}
    \put(24,0){\line(2,1){24}}
  \end{picture}
  \quad , \quad
  \begin{picture}(48,36)
    \put(0,18){\makebox(0,0){$\ast$}}
    \put(48,18){\makebox(0,0){$-$}}
     
    \put(0,12){\line(2,-1){24}}
    \put(24,0){\line(2,1){24}}
  \end{picture}
  \quad \}
  \qquad
  \entails_\typeBin
  \qquad
  \begin{picture}(48,36)
    \put(0,18){\makebox(0,0){$-$}}
    \put(48,18){\makebox(0,0){$-$}}
     
    \put(0,12){\line(2,-1){24}}
    \put(24,0){\line(2,1){24}}
  \end{picture}
\end{equation*}
and also
\begin{equation*}
  \begin{picture}(60,36)
    \put(0,18){\makebox(0,0){$\ast$}}
    \put(36,30){\makebox(0,0){$\ast$}}
    \put(60,30){\makebox(0,0){$-$}}
    
    \put(0,12){\line(2,-1){24}}
    \put(24,0){\line(2,1){24}}
    \put(48,12){\line(-1,1){12}}
    \put(48,12){\line(1,1){12}}
  \end{picture}
  \qquad
  \entails_\typeBin
  \qquad
  \begin{picture}(60,36)
    \put(0,18){\makebox(0,0){$\ast$}}
    \put(36,30){\makebox(0,0){$\ast$}}
    \put(60,30){\makebox(0,0){$\ast$}}
    
    \put(0,12){\line(2,-1){24}}
    \put(24,0){\line(2,1){24}}
    \put(48,12){\line(-1,1){12}}
    \put(48,12){\line(1,1){12}}
  \end{picture}
  \quad
  \hbox{and}
  \quad
  \begin{picture}(48,36)
    \put(0,18){\makebox(0,0){$\ast$}}
    \put(48,18){\makebox(0,0){$\ast$}}
     
    \put(0,12){\line(2,-1){24}}
    \put(24,0){\line(2,1){24}}
  \end{picture}
\end{equation*}

\begin{definition*}
  \emph{Objects} of a base type are finite or infinite sets of tokens
  with are consistent and closed under entailment.  An object $x$ is
  \emph{cototal} if for each of its tokens $P(\ast)$ with a
  distinguished occurrence of $\ast$ there is another token in $x$
  where this occurrence of $\ast$ is replaced by a constructor
  expression with only $\ast$'s as arguments.  We call $x$
  \emph{total}\index{total} if it is cototal and finite.
\end{definition*}

The tokens of $\typeN$ are shown in Figure~\ref{F:N}.  For tokens
$a,b$ we have $\{ a \} \entails b$ if and only if there is a path from
$a$ (up) to $b$ (down).  There is exactly one cototal object,
consisting of $\suc\ast$, $\suc\suc\ast$, $\suc\suc\suc\ast$ \dots.
Note that $\ast$ is not a token; the \inquotes{bottom} object is the
empty set.
\begin{figure}
  \begin{center}
    \begin{picture}(170,81)
      \put(24,0){\makebox(0,0){$\bullet$}}
      \put(12,0){\makebox(0,0){$0$}}
      \put(72,0){\makebox(0,0){$\bullet$}}
      \put(90,0){\makebox(0,0){$\suc\ast$}}

      \put(72,0){\line(-1,1){24}}
      \put(48,24){\makebox(0,0){$\bullet$}}
      \put(30,24){\makebox(0,0){$\suc 0$}}
      \put(72,0){\line(1,1){24}}
      \put(96,24){\makebox(0,0){$\bullet$}}
      \put(120,24){\makebox(0,0){$\suc \suc \ast$}}

      \put(96,24){\line(-1,1){24}}
      \put(72,48){\makebox(0,0){$\bullet$}}
      \put(45,48){\makebox(0,0){$\suc \suc 0$}}
      \put(96,24){\line(1,1){24}}
      \put(120,48){\makebox(0,0){$\bullet$}}
      \put(150,48){\makebox(0,0){$\suc \suc \suc \ast$}}

      \put(120,48){\line(-1,1){24}}
      \put(96,72){\makebox(0,0){$\bullet$}}
      \put(60,72){\makebox(0,0){$\suc \suc \suc 0$}}
      \put(120,48){\line(1,1){24}}
      \put(147,75){\makebox(0,0){.}}
      \put(150,78){\makebox(0,0){.}}
      \put(153,81){\makebox(0,0){.}}
    \end{picture}
  \end{center}
  \caption{Tokens and entailment for $\typeN$}
  \label{F:N}
\end{figure}

Examples of more interesting objects in $\typeBin$ are
\begin{enumerate}
\item $R \defeq $ closure (under entailment) of all tokens of the form
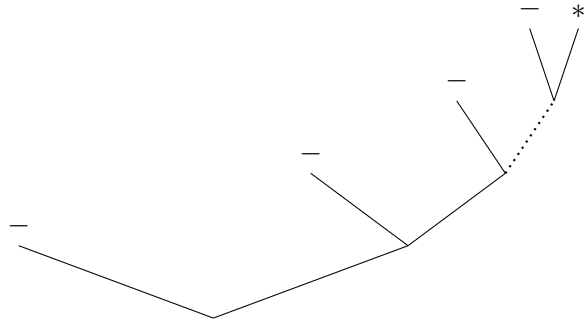
\begin{equation*}
  \begin{tikzpicture}[scale=.16]
    \coordinate (0) at (32,0);
    \coordinate [label=above:$-$] (-1/2) at (16,6);
    \coordinate (1/2) at (48,6);
    \coordinate (-3/4) at (8,12);
    \coordinate (-1/4) at (24,12);
    \coordinate [label=above:$-$] (1/4) at (40,12);
    \coordinate (3/4) at (56,12);
    \coordinate (-7/8) at (4,18);
    \coordinate (-5/8) at (12,18);
    \coordinate (-3/8) at (20,18);
    \coordinate (-1/8) at (28,18);
    \coordinate (1/8) at (36,18);
    \coordinate (3/8) at (44,18);
    \coordinate [label=above:$-$] (5/8) at (52,18);
    \coordinate (7/8) at (60,18);
    \coordinate (-15/16) at (2,24);
    \coordinate (-13/16) at (6,24);
    \coordinate (-11/16) at (10,24);
    \coordinate (-9/16) at (14,24);
    \coordinate (-7/16) at (18,24);
    \coordinate (-5/16) at (22,24);
    \coordinate (-3/16) at (26,24);
    \coordinate (-1/16) at (30,24);
    \coordinate (1/16) at (34,24);
    \coordinate (3/16) at (38,24);
    \coordinate (5/16) at (42,24);
    \coordinate (7/16) at (46,24);
    \coordinate (9/16) at (50,24);
    \coordinate (11/16) at (54,24);
    \coordinate [label=above:$-$] (13/16) at (58,24);
    \coordinate [label=above:$\ast$] (15/16) at (62,24);

    \draw (0) edge (-1/2);
    \draw (0) edge (1/2);
    \draw (1/2) edge (1/4);
    \draw (1/2) edge (3/4);
    \draw (3/4) edge (5/8);
    \draw [dotted, thick] (3/4) edge (7/8);
    \draw (7/8) edge (13/16);
    \draw (7/8) edge (15/16);
  \end{tikzpicture}
\end{equation*}

\item $L$ is defined similarly

\item $L \cup R$

\item $\frac{1}{2} \defeq$ closure of all tokens of the form
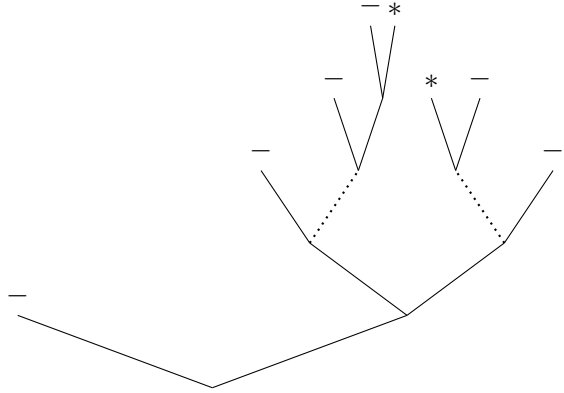
\begin{equation*}
  \begin{tikzpicture}[scale=.16]
    \coordinate (0) at (32,0);
    \coordinate [label=above:$-$] (-1/2) at (16,6);
    \coordinate (1/2) at (48,6);
    \coordinate (-3/4) at (8,12);
    \coordinate (-1/4) at (24,12);
    \coordinate (1/4) at (40,12);
    \coordinate (3/4) at (56,12);
    \coordinate (-7/8) at (4,18);
    \coordinate (-5/8) at (12,18);
    \coordinate (-3/8) at (20,18);
    \coordinate (-1/8) at (28,18);
    \coordinate [label=above:$-$] (1/8) at (36,18);
    \coordinate (3/8) at (44,18);
    \coordinate (5/8) at (52,18);
    \coordinate [label=above:$-$] (7/8) at (60,18);
    \coordinate (-15/16) at (2,24);
    \coordinate (-13/16) at (6,24);
    \coordinate (-11/16) at (10,24);
    \coordinate (-9/16) at (14,24);
    \coordinate (-7/16) at (18,24);
    \coordinate (-5/16) at (22,24);
    \coordinate (-3/16) at (26,24);
    \coordinate (-1/16) at (30,24);
    \coordinate (1/16) at (34,24);
    \coordinate (3/16) at (38,24);
    \coordinate [label=above:$-$] (5/16) at (42,24);
    \coordinate (7/16) at (46,24);
    \coordinate [label=above:$\ast$] (9/16) at (50,24);
    \coordinate [label=above:$-$] (11/16) at (54,24);
    \coordinate (13/16) at (58,24);
    \coordinate (15/16) at (62,24);

    \coordinate [label=above:$-$] (15/32) at (45,30);
    \coordinate [label=above:$\ast$] (19/32) at (47,30);
 
    \draw (0) edge (-1/2);
    \draw (0) edge (1/2);
    \draw (1/2) edge (1/4);
    \draw (1/2) edge (3/4);
    \draw (1/4) edge (1/8);
    \draw [dotted, thick] (1/4) edge (3/8);
    \draw [dotted, thick] (3/4) edge (5/8);
    \draw (3/4) edge (7/8);
    \draw (3/8) edge (5/16);
    \draw (3/8) edge (7/16);
    \draw (5/8) edge (9/16);
    \draw (5/8) edge (11/16);

    \draw (7/16) edge (15/32);
    \draw (7/16) edge (19/32);
  \end{tikzpicture}
\end{equation*}
\item Closure of
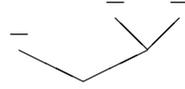
\begin{equation*}
  \begin{picture}(60,36)
    \put(0,18){\makebox(0,0){$-$}}
    \put(36,30){\makebox(0,0){$-$}}
    \put(60,30){\makebox(0,0){$-$}}
    
    \put(0,12){\line(2,-1){24}}
    \put(24,0){\line(2,1){24}}
    \put(48,12){\line(-1,1){12}}
    \put(48,12){\line(1,1){12}}
  \end{picture}
\end{equation*}
\item Closure of
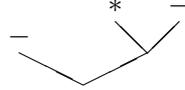
\begin{equation*}
  \begin{picture}(60,36)
    \put(0,18){\makebox(0,0){$-$}}
    \put(36,30){\makebox(0,0){$\ast$}}
    \put(60,30){\makebox(0,0){$-$}}
    
    \put(0,12){\line(2,-1){24}}
    \put(24,0){\line(2,1){24}}
    \put(48,12){\line(-1,1){12}}
    \put(48,12){\line(1,1){12}}
  \end{picture}
\end{equation*}
\end{enumerate}
Among these are
\begin{alignat*}{2}
  &(1)-(4) &\quad&\hbox{cototal objects},
  \\
  &(5) &\quad&\hbox{a total object},
  \\
  &(6) &\quad&\hbox{a general object}.
\end{alignat*}

Finally we consider the algebra $\typeL{\typeB}$ of lists of booleans,
to represent finite or infinite paths.  Tokens in $\typeL{\typeB}$ are
for instance (writing $::$ for $\mathrm{Cons}$)
\begin{equation*}
  \{ \false \} :: \{ \false \} :: \{ \true \} :: \mathrm{Nil}
  \qquad \hbox{abbreviated} \qquad
  RRL\nil
\end{equation*}
Consistency in $\typeL{\typeB}$:
\begin{equation*}
  RL\nil \not\consistent RRL\ast,
  \qquad
  RR\ast \consistent RRL\ast
\end{equation*}
Entailment in $\typeL{\typeB}$:
\begin{align*}
  &RRL\ast \entails RR\ast, R\ast
  \qquad \hbox{(and also e.g.\ $RR\emptyset \ast$)},
  \\
  &RRL\nil \entails RRL\ast, RR\ast, R\ast
\end{align*}
Objects in $\typeL{\typeB}$:
\begin{alignat*}{2}
  &\hbox{total:} &\quad& \hbox{Closure of $RLRR\nil$}
  \\
  &\hbox{cototal:} &\quad&
  \hbox{Closure of all $RLRR\dots R\ast$ and all $RRLL\dots L\ast$}
  \\
  &\hbox{general:} &\quad& \hbox{Closure of $RLRR \ast$}
\end{alignat*}

Objects of higher type are canonically defined in the Scott-Ershov
model \cite{Scott70, Ershov77}.  A detailed exposition is in
\cite[Chapter 6]{SchwichtenbergWainer12}.

\subsection{Type theory for finite and infinite data}
\label{SS:TCF}
Terms are built from (typed) variables, constructors and
(equationally) defined constants by $\lambda$-abstraction and
application, where termination of the defining equations is not
required.  Defined constants are given by their type and computation
rules.  Examples for the type $\typeList \defeq \typeL{\typeSd}$ of
lists of signed digits are

Terms are built from (typed) variables, constructors and
(equationally) defined constants by $\lambda$-abstraction and
application, where termination of the defining equations is not
required.  Defined constants are given by their type and computation
rules.  Examples for the type $\typeList \defeq \typeL{\typeSd}$ of
lists of signed digits are
\begin{alignat*}{3}
  &\rec_\typeList^\tau &&\colon \typeList \typeTo \tau \typeTo
  (\typeSd \typeTo \typeList \typeTo \tau \typeTo \tau) \typeTo \tau
  &&\qquad\hbox{(recursion operator)},
  \\
  &\destr_\typeList &&\colon \typeList \typeTo
  \typeUnit \typeSum (\typeSd \typeProd \typeList)
  &&\qquad\hbox{(destructor)},
  \\
  &\corec_\typeList^\tau &&\colon \tau \typeTo (\tau \typeTo \typeUnit
  \typeSum( \typeSd \typeProd (\typeList \typeSum \tau))) \typeTo
  \typeList &&\qquad\hbox{(corecursion operator)}.
\end{alignat*}
$\rec_\typeList^\tau$ takes a recursion argument $\ell \colon
\typeList$, an initial value $z \colon \tau$ and a \inquotes{step}
argument $f \colon (\typeSd \typeTo \typeList \typeTo \tau \typeTo
\tau)$ operating according to the structure of the given list.  The
meaning of $\rec_\typeList^\tau$ is determined by the computation
rules
\begin{equation*}
  \rec_\typeList^\tau( \nil, z, f) = z,
  \qquad
  \rec_\typeList^\tau( s:: \ell, z, f) =
  f(s, \ell, \rec_\typeList^\tau( \ell, z, f)).
\end{equation*}
$\destr_\typeList$ operates by disassembling the given list $\ell
\colon \typeList$ into its head and tail.  The computation rules are
\begin{equation*}
  \destr_\typeList(\nil) = \termSumIntroLeft(\termUnit),
  \qquad
  \destr_\typeList(s :: \ell) = \termSumIntroRight \pair {s} {\ell}.
\end{equation*}
For $\corec_\typeList^\tau$ the \inquotes{step} $f \colon (\tau
\typeTo \typeUnit \typeSum( \typeSd \typeProd (\typeList \typeSum
\tau)))$ operates by inspection (or observation) of its argument $z
\colon \tau$.  More precisely, the meaning of $\corec_\typeList^\tau$
is given by the (non-terminating) computation rule
\begin{equation*}
  \corec_\typeList^\tau z f =
  \begin{cases}
    \nil
    &\hbox{if $\eqd{fz}{\termSumIntroLeft(\termUnit)}$},
    \\
    s :: u
    &\hbox{if $\eqd{fz}{\termSumIntroRight\pair{s}{\termSumIntroLeft(u)}}$},
    \\
    s :: \corec_\typeList^\tau z' f
    &\hbox{if $\eqd{fz}{\termSumIntroRight\pair{s}{\termSumIntroRight(z')}}$}.
  \end{cases}
\end{equation*}
We use $\equiv$ for Leibniz equality (inductively defined by the
clause $\forall_z( \eqd{z}{z})$).

Predicates are either (co)inductively defined or of the form
$\set{\vec{x}}{A}$ where $A$ is a formula.  We identify
$\set{\vec{x}}{A(\vec{x}\,)} \vec{t}$ with $A(\vec{t\,}\,)$.  Formulas
are given as atomic formulas $P \vec{t}$ where $P$ is a predicate and
$\vec{t}$ are terms, or else as $A \to B$ or $\forall_x A$ where $A$
and $B$ are already formulas.  The formulas $\ex_xA$, $A \land B$ and
$A \lor B$ technically are inductively defined (nullary) predicates.

An inductively defined predicate $I$ is given by \emph{clauses} of the form
\begin{equation*}
  I_i^+ \colon \forall_{\vec{x}_i}((A_{ij}(I))_{j<n_i} \to I\vec{t}_i), 
\end{equation*}
for $i\in \{0,\dots,k-1\}$, $k>0$ and $n_i\geq 0$, where $I$ only
occurs strictly positive in each $A_{ij}(I)$ and $\vec{x_i}$ contains
all free variables of $A_{ij}(I)$ and $\vec{t}_i$.  The clauses are
called introduction axioms of $I$.

An example is the totality predicate $\Total_\typeList$.  Its clauses
are
\begin{equation*}
  (\Total_\typeList)^+_0 \colon \Total_\typeList \nil,
  \qquad
  (\Total_\typeList)^+_1 \colon \forall_{s,\ell}(
  \Total_\typeSd s \to \Total_\typeList \ell \to
  \Total_\typeList s :: \ell),
\end{equation*}
where $\Total_\typeSd$ is the totality predicate for the three-element
algebra $\typeSd$ of signed digits, defined by the three clauses
$\Total_\typeSd s$ for $s$ a signed digit.

Intuitively, an inductively defined predicate is the smallest
(w.r.t.\ $\subseteq$) predicate fulfilling its clauses.  This is
expressed by its elimination axiom:
\begin{align*}
  I^- \colon
  (\forall_{\vec{x}_i}((A_{ij}(I\cap X))_{j<n_i} \to X\vec{t}_i))_{i<k} \to
  I\subseteq X,
\end{align*}
where $X$ can be any predicate with the same arity as $I$.  The axiom
$I^-$ is also called induction axiom.  For example,
\begin{equation*}
  \Total_\typeList^- \colon X \nil \to \forall_{s,\ell}(
    \Total_\typeSd s \to
    \Total_\typeList \ell \to X \ell \to 
    X(s :: \ell)) \to \Total_\typeList \subseteq X.
\end{equation*}

Let $I$ be an inductively defined predicate given by its clauses
$(I^+_i)_{i<k}$, and assume that there are recursive calls in some
clauses.  To every such $I$ we associate a \inquotes{companion} $\CoI$
and call it a coinductively defined predicate.  The elimination (or
closure) axiom $\CoI^-$ of the coinductively defined predicate says
that if $\CoI\vec{x}$ holds, then it comes from one of the clauses:
\begin{equation*}
  \CoI^- \colon 
  \forall_{\vec{x}}(\CoI \vec{x} \to \biglor_{i<k} \ex_{\vec{x}_i}(
  \bigland_{j < n_i} A_{i j}(\CoI) \land \eqd{\vec{x}} {\vec{t}_i})).
\end{equation*}
In words: if $\CoI \vec{x}$ holds, then there is at least one clause
$I^+_i$ whose premises $(A_{ij}(\CoI))_{j<n_i}$ (with $\CoI$ instead
of $I$) are fulfilled and whose conclusion is $\CoI\vec{x}$ up to
Leibniz equality.  For example, the co-predicate associated with
$\Total_\typeList$ is cototality $\CoTotal_\typeList$.  Its
elimination axiom is
\begin{equation*}
  \CoTotal_\typeList^- \colon \CoTotal_\typeList \ell \to
  \ell \equiv \nil \lor \ex_{s,\ell'}(
  \Total_\typeSd s \land \CoTotal_\typeList \ell' \land
  \ell \equiv s :: \ell').
\end{equation*}

The introduction (or coinduction or greatest-fixed-point) axiom
$\CoI^+$ of $\CoI$ says that $\CoI$ is the greatest predicate such that
\begin{equation*}
  \CoI^+ \colon \forall_{\vec{x}}( X \vec{x} \to
  \biglor_{i<k} \ex_{\vec{x}_i}(
  \bigland_{j < n_i} A_{i j}(\CoI \cup X) \land
  \eqd{\vec{x}} {\vec{t}_i})) \to X \subseteq \CoI.
\end{equation*}
This means that any other \inquotes{competitor} predicate $X$
satisfying the same closure property is below $\CoI$.  For example,
\begin{equation*}
  \CoTotal_\typeList^+ \colon \forall_\ell(X \ell \to
  \ell \equiv \nil \lor \ex_{s,\ell'}(
  \Total_\typeSd s \land
  (\CoTotal_\typeList \ell' \lor X \ell') \land
  \ell \equiv s :: \ell') \to X \subseteq {\CoTotal_\typeList}.
\end{equation*}

\subsection{Computational content}
\label{SS:CompContent}
We are interested in the computational content of proofs, which in the
present setting arises from inductively and coinductively defined
predicates only: they can be declared to be computationally relevant
(c.r.)\ or non-computational (n.c.).  This also applies to predicate
variables $X$.  The (finitely many) clauses of a c.r.\ inductive
predicate $I$ determine an algebra $\iota_I$.  A \inquotes{witness}
that $I \vec{t}$ or $\CoI \vec{t}$ holds has type $\iota_I$.  It will
be a total object in the first and a cototal object in the second
case.  For instance, the type of the c.r.\ version of $\Total_\typeList$
or $\CoTotal_\typeList$ is $\typeList$.  A formula $A$ built from
atomic formulas by $\to$, $\forall$ is c.r.\ if its final conclusion
is.  The \emph{type} $\extrTy{A}$ is defined by
\begin{align*}
  \extrTy{P \vec{t}\;} &\defeq \extrTy{P},
  \\
  \extrTy{A \to B} &\defeq
  \begin{cases}
    \extrTy{A} \typeTo \extrTy{B} & \hbox{if $A$ is c.r.}
    \\
    \extrTy{B} & \hbox{if $A$ is n.c.}
  \end{cases}
  \\
  \extrTy{\forall_x A} &\defeq \extrTy{A}.
\end{align*}
Recall that $\ex_xA$, $A \land B$ and $A \lor B$ technically are
inductively defined (nullary) predicates.  If they are c.r.,
their types are $\extrTy{\ex_x A} \defeq \extrTy{A}$; $\extrTy{A \land
  B} \defeq \extrTy{A} \typeProd \extrTy{B}, \extrTy{A}, \extrTy{B}$
depending on the c.r./n.c.\ status of $A,B$; $\extrTy{A \lor
  B} \defeq \extrTy{A} \typeSum \extrTy{B}, \extrTy{A} \typeSum
\typeUnit, \typeUnit \typeSum \extrTy{B}, \typeB$ depending on the
c.r./n.c.\ status of $A,B$.

Let $M$ be a proof in $\TCF$ of a c.r.\ formula $A$.  We define its
\emph{extracted term} $\extrTer{M}$, of type $\extrTy{A}$, with the
aim to express $M$'s computational content.  It will be a term built
up from variables, constructors, recursion operators, destructors
and corecursion operators by $\lambda$-abstraction and application:
\begin{alignat*}{2}
  &\extrTer{u^A} &&\defeq z_u^{\extrTy{A}} \quad\hbox{
    ($z_u^{\extrTy{A}}$ uniquely associated with $u^A$)},
  \\
  &\extrTer{(\lambda_{u^A} M^B)^{A \to B}} &&\defeq
  \begin{cases}
    \lambda_{z_u} \extrTer{M} &\hbox{if $A$ is c.r.}
    \\
    \extrTer{M} &\hbox{if $A$ is n.c.}
  \end{cases}
  \\
  &\extrTer{(M^{A \to B}N^A)^B} &&\defeq
  \begin{cases}
    \extrTer{M}\extrTer{N} &\hbox{if $A$ is c.r.}
    \\
    \extrTer{M} &\hbox{if $A$ is n.c.}
  \end{cases}
  \\
  &\extrTer{(\lambda_x M^A)^{\forall_x A}} &&\defeq \extrTer{M},
  \\
  &\extrTer{(M^{\forall_x A(x)} t)^{A(t)}} &&\defeq \extrTer{M}.
\end{alignat*}
It remains to define extracted terms for the axioms.  Consider a
(c.r.)\ inductively defined predicate $I$.  For its introduction and
elimination axioms define $\extrTer{I^+_i} \defeq \constr_i$ and
$\extrTer{I^-} \defeq \rec$, where both the constructor $\constr_i$
and the recursion operator $\rec$ refer to the algebra $\iota_I$
associated with $I$.  For the closure and greatest-fixed-point axioms
of $\CoI$ define $\extrTer{\CoI^-} \defeq \destr$ and
$\extrTer{\CoI^+} \defeq \corec$, where both the destructor $\destr$
and the corecursion operator $\corec$ again refer to the algebra
$\iota_I$ associated with $I$.

\subsection{Real numbers}
\label{SS:Reals}
Rational numbers $\typeQ$ are defined by $\typeQ \defeq
\typeZ\times \typeP$, where the pair of $k \colon \typeZ$ and $p
\colon \typeP$ is denoted by $\frac{k}{p}$.  We also have the expected
notion of equality between rational numbers.  In contrast to the
positive, natural, rational numbers or the integers, the real numbers
are not an algebra.  The property to be a real number is a predicate
on $(\typeN \typeTo \typeQ) \typeProd (\typeP \typeTo \typeN)$.
Totality on this type means that both components map total arguments
into total values.  By a real number $x$ we mean a pair of this type
whose first component is a Cauchy sequence $(a_n)_{n \in \typeN}$ of
rationals with the second component $M$ as modulus, i.e.,
\begin{equation*}
  |a_n - a_m| \le \frac{1}{2^p} \quad \hbox{for $n,m \ge M(p)$}.
\end{equation*}
We write $\typeR$ for this predicate.  For real numbers we
can define the arithmetic functions $+$, $-$, $\cdot$ and the inverse
and also the absolute value $|\cdot|$.  Note that the inverse needs a
witness that the real number is positive.  We will not go into details
here.  An exact definition of these functions and also of positivity
can be found in \cite{Wiesnet17, wiesnet18} and in the library of
Minlog\footnote{\url{http://minlog-system.de}.}.

An important predicate on the real numbers is real equality: we define
$x \le y$ by
\begin{equation*}
  a_{M(p+1)} \le b_{N(p+1)} + \frac{1}{2^p} \quad
  \hbox{for all $p \in \typeP$}
\end{equation*}
and real equality $x=y$ by $x \le y \land y \le x$.  Most of the
predicates and functions on real numbers will be compatible with real
equality; however, this must be proved in each case.

The rational numbers can be embedded into the real numbers by
identifying a rational number $a$ with the pair consisting of the
constant Cauchy sequence $(a)_{n \in \typeN}$ and the constant modulus
$(0)_{p\in \typeP}$.

We use variable names
\begin{alignat*}{2}
  &x,y
  &\quad&\hbox{for real numbers},
  \\
  &a,b,c
  &\quad&\hbox{for rational numbers},
  \\
  &d,e,k,j,i
  &\quad&\hbox{for integers},
  \\
  &s
  &\quad&\hbox{for signed digits},
  \\
  &u,v
  &\quad&\hbox{for lists (and streams) of signed digits}.
\end{alignat*}
Let $\ncind{\Total_\typeZ}$ be the n.c.\ version of the totality
predicate for $\typeZ$, and similarly for other types.  Let
$\mathrm{Sd}$ be the (formally inductive) predicate consisting of the
integers $1,0,-1$, i.e., $\mathrm{Sd}$ is given by the three clauses
$1\in \mathrm{Sd}$, $0\in \mathrm{Sd}$ and $-1\in \mathrm{Sd}$.  We
require that $\mathrm{Sd}$ has computational content (in the
three-element algebra $\typeSd$), and write
\begin{equation*}
  \hbox{$\forall_{d \in \mathrm{Sd}} A(d)$ \quad for \quad
    $\forall_{\hat{d}}( \hat{d} \in \ncind{\Total_{\mathbb{Z}}} \to
    \hat{d} \in \mathrm{Sd} \to A(\hat{d}))$}
\end{equation*}
and similarly for $\ex$.  In this case the premise $\hat{d} \in
\ncind{\Total_\typeZ}$ is not necessary because it follows from $\hat{d}
\in \mathrm{Sd}$.  The predicate $\typeR$ and also $\le, =$ on real
numbers are assumed to be n.c.  We write (with $x$ ranging over
$[-1,1]$)
\begin{align*}
  \hbox{$\forall_x A(x)$} &\quad \hbox{for} \quad
    \hbox{$\forall_{\hat{x}}(
      \hat{x} \in \typeR \to
      A(\hat{x}))$}
  \\
  \hbox{$\forall_{x \in P} A(x)$} &\quad \hbox{for} \quad
    \hbox{$\forall_{\hat{x}}( 
      \hat{x} \in \typeR \to
      \hat{x} \in P \to A(\hat{x}))$}
\end{align*}
and again similarly for $\ex$.  The totality of $\hat{x}$ is part of
$\hat{x} \in \typeR$.

We inductively define a predicate $I_0$ on reals $s$ such that the
generating sequence $d_0, \dots, d_{m-1}$ for $x \in I_0$ means that
$x$ is between $\sum_{n<m} \frac{d_n}{2^{n+1}} - \frac{1}{2^m}$ and
$\sum_{n<m} \frac{d_n}{2^{n+1}} + \frac{1}{2^m}$.  Let $I_0$ be
defined by the two clauses
\begin{equation*}
  \forall_x(x=0 \to x \in I_0), \quad
  \forall_{d \in \mathrm{Sd}} \forall_x\forall_{x' \in I_0}\Bigl(
    x=\frac{x'+d}{2}\to x \in I_0 \Bigr).
\end{equation*}
Then the induction axiom is
\begin{equation*}
  \forall_x(x=0\to x\in P) \to
  \forall_{d \in \mathrm{Sd}} \forall_x\forall_{x' \in I_0 \cap P}\Bigl(
  x= \frac{x+d}{2}\to x \in P \Bigr) \to I_0 \subseteq P.
\end{equation*}

The clauses of $I_0$ also determine the dual predicate $\CoI_0$, which
is coinductively given by the closure axiom
\begin{equation*}
  \forall_{x \in \CoI_0} \Bigl(
  x=0 \lor \ex_{d \in \mathrm{Sd}} \ex_{x' \in \CoI_0}
  \Bigl( x=\frac{x'+d}{2} \Bigr) \Bigr)
\end{equation*}
and the coinduction axiom 
\begin{equation*}
  \forall_{x \in P} \Bigl( x=0 \lor 
  \ex_{d \in \mathrm{Sd}} \ex_{x' \in \CoI_0 \cup P}
  \Bigl( x=\frac{x'+d}{2} \Bigr) \Bigr) \to P \subseteq \CoI_0.
\end{equation*}

We require that both predicates $I_0$ and $\CoI_0$ have computational
content.  The data type associated to these clauses is the algebra
$\typeList$ of lists of signed digits.  This association is quite
canonical: the first constructor $\nil \colon \typeList$ is a witness
for the first clause, which does not need any further information. The
second constructor $::$ of type $\typeSd \typeTo \typeList \typeTo
\typeList$ is a witness for the second clause consisting of a signed
digit and a witness that $x' \in I_0$.  We refer to \cite[Section
  7.2]{SchwichtenbergWainer12} for a formal definition of the type of
a (co)inductively defined predicate.

The computational content of the clauses are the constructors, of the
induction axiom the recursion operator $\rec_\typeList^\tau$, of the
closure axiom the destructor $\destr_\typeList$ and of the coinduction
axiom the corecursion operator $\corec_\typeList^\tau$.  All these
have been explained in Section~\ref{SS:TCF}.

\section{Division for stream-represented real numbers}
\label{S:Div}

\subsection{Division for signed digit streams}
\label{SS:DivSd}
From a computational point of view, a real number $x$ is best seen as
a device producing for a given accuracy $\frac{1}{2^n}$ a rational
number $a_n$ being $\frac{1}{2^n}$-close to $x$, i.e., $|a_n - x| \le
\frac{1}{2^n}$.  For simplicity we again restrict ourselves to real
numbers in the interval $[-1,1]$ (rather than $[-2^n, 2^n]$), and to
dyadic rationals $\sum_{n<m} \frac {d_n} {2^{n+1}}$ ($d_n \in
\{ 1, \bar{1}\}$).

\begin{figure}[ht]
  \begin{tikzpicture}[scale=.16]
    \coordinate [label=below:\textcolor{red}{$0$}] (0) at (32,0);
    \coordinate [label=below:\textcolor{red}{$-\frac{1}{2}$}] (-1/2) at (16,6);
    \coordinate [label=below:\textcolor{red}{$\frac{1}{2}$}] (1/2) at (48,6);
    \coordinate [label=below:\textcolor{red}{$-\frac{3}{4}$}] (-3/4) at (8,12);
    \coordinate (-1/4) at (24,12);
    \coordinate (1/4) at (40,12);
    \coordinate [label=below:\textcolor{red}{$\frac{3}{4}$}] (3/4) at (56,12);
    \coordinate [label=left:\textcolor{red}{$-\frac{7}{8}$}] (-7/8) at (4,18);
    \coordinate (-5/8) at (12,18);
    \coordinate (-3/8) at (20,18);
    \coordinate (-1/8) at (28,18);
    \coordinate (1/8) at (36,18);
    \coordinate (3/8) at (44,18);
    \coordinate (5/8) at (52,18);
    \coordinate [label=right:\textcolor{red}{$\frac{7}{8}$}] (7/8) at (60,18);
    \coordinate[label=left:\textcolor{red}{$-\frac{15}{16}$}](-15/16) at (2,24);
    \coordinate (-13/16) at (6,24);
    \coordinate (-11/16) at (10,24);
    \coordinate (-9/16) at (14,24);
    \coordinate (-7/16) at (18,24);
    \coordinate (-5/16) at (22,24);
    \coordinate (-3/16) at (26,24);
    \coordinate (-1/16) at (30,24);
    \coordinate (1/16) at (34,24);
    \coordinate (3/16) at (38,24);
    \coordinate (5/16) at (42,24);
    \coordinate (7/16) at (46,24);
    \coordinate (9/16) at (50,24);
    \coordinate (11/16) at (54,24);
    \coordinate (13/16) at (58,24);
    \coordinate[label=right:\textcolor{red}{$\frac{15}{16}$}](15/16) at (62,24);
    
    \draw (0) edge node [below] {$\bar{1}$} (-1/2);
    \draw (0) edge node [below] {$1$} (1/2);
    \draw (-1/2) edge node [below] {$\bar{1}$} (-3/4);
    \draw (-1/2) edge node [below] {$1$} (-1/4);
    \draw (1/2) edge node [below] {$\bar{1}$} (1/4);
    \draw (1/2) edge node [below] {$1$} (3/4);
    \draw (-3/4) edge node [left] {$\bar{1}$} (-7/8);
    \draw (-3/4) edge node [right] {$1$} (-5/8);
    \draw (-1/4) edge node [left] {$\bar{1}$} (-3/8);
    \draw (-1/4) edge node [right] {$1$} (-1/8);
    \draw (1/4) edge node [left] {$\bar{1}$} (1/8);
    \draw (1/4) edge node [right] {$1$} (3/8);
    \draw (3/4) edge node [left] {$\bar{1}$} (5/8);
    \draw (3/4) edge node [right] {$1$} (7/8);
    \draw (-7/8) edge node [left] {$\bar{1}$} (-15/16);
    \draw (-7/8) edge node [right] {$1$} (-13/16);
    \draw (-5/8) edge node [left] {$\bar{1}$} (-11/16);
    \draw (-5/8) edge node [right] {$1$} (-9/16);
    \draw (-3/8) edge node [left] {$\bar{1}$} (-7/16);
    \draw (-3/8) edge node [right] {$1$} (-5/16);
    \draw (-1/8) edge node [left] {$\bar{1}$} (-3/16);
    \draw (-1/8) edge node [right] {$1$} (-1/16);
    \draw (1/8) edge node [left] {$\bar{1}$} (1/16);
    \draw (1/8) edge node [right] {$1$} (3/16);
    \draw (3/8) edge node [left] {$\bar{1}$} (5/16);
    \draw (3/8) edge node [right] {$1$} (7/16);
    \draw (5/8) edge node [left] {$\bar{1}$} (9/16);
    \draw (5/8) edge node [right] {$1$} (11/16);
    \draw (7/8) edge node [left] {$\bar{1}$} (13/16);
    \draw (7/8) edge node [right] {$1$} (15/16);
  \end{tikzpicture}
  \caption{Dyadic rationals.}
\end{figure}

Clearly we can represent real numbers (in $[-1,1]$) as streams of
$1, \bar{1}$.  Now when doing arithmetic operations on such streams
the problem of \inquotes{productivity} arises: suppose we want to add
the two streams $\bar{1} 1 1 1 \dots$ and $1 \bar{1} \bar{1} \bar{1}
\dots$.  Then the first digit of the output stream will only be known
after we have checked the two input streams long enough, but there is
no bound how far we have to look.  The well-known cure for this
problem is to add a \emph{delay} digit $0$ and work with signed digits
$\sum_{n<m} \frac {d_n} {2^{n+1}}$, now with $d_n \in
\{1,0,\bar{1}\}$.  We have a lot of redundancy here (for instance
$\bar{1}1$ and $0 \bar{1}$ both denote $-\frac{1}{4}$), but this is
not a serious problem.

Since $0$ as real number is represented by the stream consisting of
the digit $0$ only, we can simplify our example by removing the
nullary clause from the inductive definition of $I_0$, and define $I$
and $\CoI$ accordingly.  We only need $\CoI$, coinductively defined by
the closure axiom
\begin{equation}
  \label{E:ClosureI}
  \forall_{x \in \CoI} \ex_{d \in \mathrm{Sd}} \ex_{x' \in \CoI}
  \Bigl( x=\frac{x'+d}{2} \Bigr).
\end{equation}
Therefore, the coinduction axiom is 
\begin{equation}
  \label{E:GfpI}
  \forall_{x \in P} \ex_{d \in \mathrm{Sd}} \ex_{x' \in \CoI \cup P}
  \Bigl( x=\frac{x'+d}{2} \Bigr) \to P \subseteq \CoI.
\end{equation}
The canonically associated data type then becomes the algebra
$\typeSds$ given by a single binary constructor of type $\typeSd
\typeTo \typeSds \typeTo \typeSds$, again denoted by $::$ (infix).
Note that we do not require that $\Total_\typeIntv u$ holds for any
$u$.  Objects $u \in \CoTotal_\typeIntv$ are called \emph{streams} of
signed digits.

The computational content of the closure axiom \eqref{E:ClosureI} now
is the destructor $\destr_\typeSds$ of type $\typeSds \typeTo \typeSd
\typeProd \typeSds$, defined by
\begin{equation*}
  \destr_\typeSds(s :: u) = \pair {s} {u}.
\end{equation*}
For the coinduction axiom \eqref{E:GfpI} we have the corecursion
operator $\corec_\typeSds^\tau$ of type $\tau \typeTo (\tau \typeTo
\typeSd \typeProd (\typeSds \typeSum \tau)) \typeTo \typeSds$.  Note
that $\typeSd \typeProd (\typeSds \typeSum \tau)$ appears since
$\typeSds$ has the a single constructor of type $\typeSd \typeTo
\typeSds \typeTo \typeSds$.  The meaning of $\corec_\typeSds^\tau$ is
determined by a (non-terminating) computation rule similar to the one
for $\corec_\typeList^\tau$ (see Section~\ref{SS:TCF}), with the first
case left out.

Our goal in this section is to prove that $[-1,1]$ is closed under
division under some conditions, w.r.t.\ the representation of reals as
signed digit streams.  For division $\frac{x}{y}$ we clearly need a
restriction on the denominator $y$ to stay in the interval $[-1,1]$,
and must assume that $|y|$ is strictly positive.  For simplicity we
assume $\frac{1}{4} \le y$; this can easily be extended to the case
$2^{-p} \le y$ (using induction on $p$ and Lemma~\ref{L:Double}
below).  The main idea of the algorithm and also of the proof, which
is inspired by \cite{CiaffaglioneGianantonio06}, consists in three
representations of $\frac{x}{y}$:
\begin{equation*}
  \frac{x}{y} = \frac{1+\frac{x_1}{y}}{2}
  = \frac{0+\frac{x_0}{y}}{2}
  = \frac{-1+\frac{x_{-1}}{y}}{2}
\end{equation*}
where
\begin{equation*}
  x_1 = 4\frac{x+\frac{-y}{2}}{2},\quad
  x_0 = 2x,\quad
  x_{-1} = 4\frac{x+\frac{y}{2}}{2}.
\end{equation*}
Depending on what we know about $x$ we choose one of these
representations of $\frac{x}{y}$ to obtain its first digit.  This will
give us a corecursive definition of $\frac{x}{y}$.

As we see in the representation of $\frac{x}{y}$, we will use that
$\CoI$ is closed under the average:

\begin{thm}
  \label{T:CoIAverage}
  The average of two real numbers $x,y$ in $\CoI$ is in $\CoI$:
  \begin{equation*}
    \forall_{x,y \in \CoI}\Bigl( \frac{x+y}{2} \in \CoI \Bigr).
  \end{equation*}
\end{thm}

\begin{proof}
  The proof in \cite{Berger09,BergerMiyamotoSchwichtenbergTsuiki16}
  can be adapted to the present setting with concrete rather than
  abstract reals.
\end{proof}

The term extracted from this proof corresponds to an algorithm
transforming stream representations of $x$ and $y$ into a stream
representation of their average $\frac{x+y}{2}$.  For the first $n$
digits in the representation of $\frac{x+y}{2}$ one needs the first
$n+1$ digits in the representations of $x$ and $y$.

\begin{lem}
  \label{L:Shift}
  $\CoI$ is closed under shifting a real $x \le 0$ ($x \ge 0$) by $+1$
  ($-1$):
  \begin{align*}
    &\forall_{x \in \CoI}( x \le 0 \to x+1 \in \CoI),
    \\
    &\forall_{x \in \CoI}( 0 \le x \to x-1 \in \CoI).
  \end{align*}
\end{lem}

\begin{proof}
  We only consider the first claim; the second one is proved
  similarly.  Since we want to prove that something is in $\CoI$, the
  proof must use its introduction axiom given in \eqref{E:GfpI}:
  \begin{equation*}
    \forall_{x \in P} \ex_{d \in \mathrm{Sd}} \ex_{x' \in \CoI \cup P}
    \Bigl( x=\frac{x'+d}{2} \Bigr) \to \forall_x(x\in P \to x\in \CoI).
  \end{equation*}
  The crucial point is how to choose the competitor predicate $P$.
  Looking at the formula we want to prove, we take $P \defeq
  \set{x}{\ex_{y \in \CoI}(y \le 0 \land x=y+1)}$ as competitor
  predicate.  Then we have
  \begin{equation*}
    \forall_{x \in P} \ex_{d \in \mathrm{Sd}} \ex_{x' \in \CoI \cup P}
    \Bigl( x=\frac{x'+d}{2} \Bigr) \to
    \forall_x(\ex_{y \in \CoI}(y \le 0 \land x=y+1) \to x\in \CoI)
  \end{equation*}
  and since $\CoI$ is compatible with the real equality, this implies
  \begin{equation*}
    \forall_{x \in P} \ex_{d \in \mathrm{Sd}} \ex_{x' \in \CoI \cup P}
    \Bigl( x=\frac{x'+d}{2} \Bigr) \to \forall_{x\in\CoI}(x \le 0
    \to x+1\in \CoI).
  \end{equation*}
  Therefore it suffices to prove the premise of the formula above.
  Let $x\in P$ be given.  Then by definition of $P$ there is $y$ with
  $y \in \CoI$, $y\le 0$ and $x=y+1$.  From $y \in \CoI$ we know $|y|
  \le 1$ and with $y \le 0$ also $|x| \le 1$.  We need $d \in
  \mathrm{Sd}$ and $x'$ such that
  \begin{equation*}
    (x' \in \CoI \lor \ex_{y \in \CoI}(y \le 0 \land x'=y+1)) \land
    x=\frac{d+x'}{2}.
  \end{equation*}
  Again from $y \in \CoI$ we obtain $e \in \mathrm{Sd}$ and $z \in
  \CoI$ with $y=\frac{e+z}{2}$.  We now distinguish cases on $e \in
  \mathrm{Sd}$.

  \emph{Case} $e=1$.  Then $0 \ge y = \frac{1+z}{2} \ge \frac{1}{2} -
  \frac{1}{2} = 0$ and hence $x = y+1 = 1$.  Picking $d=1$ and $x'=1$
  gives the claim.  Here we need a lemma \texttt{CoIOne} stating that
  $1$ is in $\CoI$.  The proof of this lemma (by coinduction) is
  omitted.

  \emph{Case} $e=0$.  Pick $d=1$ and $x'=z+1$.  Then $z = 2y \le 0$
  and hence the r.h.s.\ of the disjunction above holds with $z$ for
  $y$.  Also $x = 2y = \frac{1+2y+1}{2}$.

  \emph{Case} $e=-1$.  Pick $d=1$ and $x'=z$.  Then the l.h.s.\ of
  the disjunction above holds, and $\frac{d+x'}{2} = \frac{1+z}{2} =
  \frac{-1+z}{2}+1 = y+1 = x$.
\end{proof}

We have formalized this proof in the Minlog\footnote{The development
  described here resides in \texttt{minlog/examples/analysis}, file
  \texttt{sddiv.scm}} proof assistant.  Minlog has a tool to extract
from the proof a term representing its computational content; it is a
term in the theory $\TCF$ described in Section~\ref{SS:TCF}.  The
extracted term is displayed as
\begin{verbatim}
[u](CoRec ai=>ai)u
 ([u0][case (DesYprod u0)
     (s pair u1 -> [case s
       (SdR -> SdR pair InL cCoIOne)
       (SdM -> SdR pair InR u1)
       (SdL -> SdR pair InL u1)])])
\end{verbatim}
Here \texttt{[u]} means lambda abstraction $\lambda_u$, and
\texttt{(CoRec ai=>ai)} is the corecursion operator
$\corec_\typeIntv^\tau$ defined above, where $\tau$ is $\typeIntv$
again.  The type of $\corec_\typeIntv^\typeIntv$ is $\typeIntv \typeTo
(\typeIntv \typeTo \typeSd \typeProd (\typeIntv \typeSum \typeIntv))
\typeTo \typeIntv$.  The first argument of the corecursion operator is
the abstracted variable \texttt{u} of type $\typeIntv$, and the second
\texttt{([u0][case \dots])} is the \inquotes{step} function, which
first destructs its argument \texttt{u0} into a pair of a signed digit
\texttt{s} and another stream \texttt{u1}, and then distinguishes
cases on \texttt{s}.  In the \texttt{SdR} case the returned pair of
type $\typeSd \typeProd (\typeIntv \typeSum \typeIntv)$ has
\texttt{SdR} again as its left component, and as right component the
$\termSumIntroLeft$ (left embedding into a sum type) of a certain
stream denoted \texttt{cCoIOne}.  This is the computational content
(hence the \inquotes{\texttt{c}}) of \texttt{CoIOne},
essentially an infinite sequence of the digit $\mathrm{SdR}$ (written
$\vec{1}$\,).

The algorithm represented by the extracted term can be understood as
follows.  If $s= \mathrm{SdR}$, then $y$ must be non-negative.  Hence
$y=0$, and a stream-representation of $y+1$ is $\vec{1}$.  Here we do
not need a corecursive call, and hence the result is $\pair
{\mathrm{SdR}} {\termSumIntroLeft(\vec{1}\;)}$.  The same happens in
case $s= \mathrm{SdL}$.  Here $y+1$ can be determined easily by
changing the first digit from $\mathrm{SdL}$ to $\mathrm{SdR}$ and
leaving the tail as it is.  Hence the result is $\pair {\mathrm{SdR}}
{\termSumIntroLeft(u')}$.  Only in case $s=\mathrm{SdM}$ we have a
corecursive call.  Here we change the first digit from $\mathrm{SdM}$
to $\mathrm{SdR}$, which however amounts to adding $\frac{1}{2}$ only.
Therefore the procedure continues with the tail.  Hence the result is
$\pair {\mathrm{SdR}} {\termSumIntroRight(u')}$.  Using the
computation rule for $\corec_\typeSds^\typeSds$ we can now describe
the computational content as a function $\texttt{add1} \colon \typeSds
\typeTo \typeSds$ defined by
\begin{align*}
  \mathrm{add1}(\mathrm{SdR}::u) &\defeq [\mathrm{SdR},\mathrm{SdR},\dots],
  \\
  \mathrm{add1}(\mathrm{SdM}::u) &\defeq \mathrm{SdR}::\mathrm{add1}(u),
  \\
  \mathrm{add1}(\mathrm{SdL}::u) &\defeq \mathrm{SdR}::u.
\end{align*}
A similar argument for the second part of the lemma gives
$\mathrm{sub1} \colon \typeSds \typeTo \typeSds$ with
\begin{align*}
  \mathrm{sub1}(\mathrm{SdR}::u) &\defeq \mathrm{SdL}::u,
  \\
  \mathrm{sub1}(\mathrm{SdM}::u) &\defeq \mathrm{SdL}::\mathrm{sub1}(u),
  \\
  \mathrm{sub1}(\mathrm{SdL}::u) &\defeq [\mathrm{SdL},\mathrm{SdL},\dots].
\end{align*}

\begin{lem}
  \label{L:Double}
  For $x$ in $\CoI$ with $|x| \le \frac{1}{2}$ we have $2x$ in $\CoI$:
  \begin{equation*}
    \forall_{x \in \CoI}\Bigl( |x| \le \frac{1}{2} \to 2x \in \CoI \Bigr).
  \end{equation*}
\end{lem}

\begin{proof}
  Let $x \in \CoI$ be given.  From the closure axiom
  \eqref{E:ClosureI} for $\CoI$ we obtain $d \in \mathrm{Sd}$ and
  $x' \in \CoI$ such that $x=\frac{d+x'}{2}$.  We distinguish cases on
  $d \in \mathrm{Sd}$.

  \emph{Case} $d=1$.  Then $x=\frac{1+x'}{2}$ and hence $2x=1+x'$.  Since
  $|x| \le \frac{1}{2}$ we have $x' \le 0$.  Now the first part of
  Lemma~\ref{L:Shift} gives the claim.

  \emph{Case} $d=0$.  Then $x=\frac{0+x'}{2}$ and hence $2x=x' \in \CoI$.

  \emph{Case} $d=-1$.  Then $x=\frac{-1+x'}{2}$ and hence $2x=-1+x'$.
  Since $|x| \le \frac{1}{2}$ we have $0 \le x'$.  Now the second part
  of Lemma~\ref{L:Shift} gives the claim.
\end{proof}

Using arguments similar to those in the remark after
Lemma~\ref{L:Shift} we can see that the corresponding algorithm can be
written as a function $\mathrm{Double} \colon \typeSds \typeTo
\typeSds$ with
\begin{align*}
  \mathrm{Double}(\mathrm{SdR}::u) &\defeq \mathrm{add1}(u),
  \\
  \mathrm{Double}(\mathrm{SdM}::u) &\defeq u,
  \\
  \mathrm{Double}(\mathrm{SdL}::u) &\defeq \mathrm{sub1}(u)
\end{align*}
where $\mathrm{add1}$ and $\mathrm{sub1}$ are the functions from this
remark.

\begin{lem}
  \label{L:CoIAux}
  For $x,y$ in $\CoI$ with $\frac{1}{4} \le y$, $|x| \le y$ and
  $0 \le x$ ($x \le 0$) we have $2x-y$ (\,$2x+y$) in $\CoI$:
  \begin{align*}
    &\forall_{x,y \in \CoI} \Bigl(
    \frac{1}{4} \le y \to |x| \le y \to 0\le x \to
    4\frac{x+\frac{-y}{2}}{2} \in \CoI \Bigr),
    \\
    &\forall_{x,y \in \CoI} \Bigl(
    \frac{1}{4} \le y \to |x| \le y \to x\le 0 \to
    4\frac{x+\frac{y}{2}}{2} \in \CoI \Bigr).
  \end{align*}
\end{lem}

\begin{proof}
  In the formulas above instead of $2x\pm y$ we have written
  $4\frac{x+(\pm y/2)}{2}$ to make Theorem~\ref{T:CoIAverage}
  applicable.  We also use two lemmas stating that $\CoI$ is closed
  under $x \mapsto \frac{x}{2}$ and $x \mapsto -x$.  To prove the
  first claim, let $x,y$ in $\CoI$ with $\frac{1}{4} \le y$, $|x| \le
  y$ and $0 \le x$.  Then clearly $\frac{-y}{2} \in \CoI$, and
  $\frac{x+\frac{-y}{2}}{2} \in \CoI$ by Theorem~\ref{T:CoIAverage}.
  We have
  \begin{equation}
    \label{E:ShiftIneq}
    \Bigl|\frac{x+\frac{-y}{2}}{2}\Bigr| = \Bigr| \frac{2x-y}{4} \Bigr|
    \le \Bigr| \frac{2y-y}{4} \Bigr| = \Bigr| \frac{y}{4} \Bigr| \le
    \frac{1}{4}.
  \end{equation}
  Hence we can apply Lemma~\ref{L:Double} twice and obtain
  $4\frac{x+\frac{-y}{2}}{2} \in \CoI$.  The proof of the second claim
  is similar.
\end{proof}

The computational content of the proofs is $\mathrm{AuxL},
\mathrm{AuxR} \colon \typeSds \typeTo \typeSds \typeTo \typeSds$:
\begin{align*}
  \mathrm{AuxL}(u,v) &\defeq
  \mathrm{Double}(\mathrm{Double}( \mathrm{Av}(u, h(-v)))),
  \\
  \mathrm{AuxR}(u,v) &\defeq
  \mathrm{Double}(\mathrm{Double}( \mathrm{Av}(u, h(v)))).
\end{align*}
Here $h,- \colon \typeSds \typeTo \typeSds$ represent the
computational content of the two lemmas used in the proof; both are
proved by coinduction.  The function $h$ prepends $\mathrm{SdM}$ to
the stream, and $-$ negates all digits:
\begin{equation*}
  \begin{split}
    h(u) \defeq \mathrm{SdM}::u,
  \end{split}
  \qquad
  \begin{split}
    -(\mathrm{SdR}::u) &\defeq \mathrm{SdL}::(-u),
    \\
    -(\mathrm{SdM}::u) &\defeq \mathrm{SdM}::(-u),
    \\
    -(\mathrm{SdL}::u) &\defeq \mathrm{SdR}::(-u).
  \end{split}
\end{equation*}
  
\begin{thm}
  \label{T:CoIDiv}
  For $x,y$ in $\CoI$ with $\frac{1}{4} \le y$ and $|x| \le y$ we have
  $\frac{x}{y}$ in $\CoI$:
  \begin{equation*}
    \forall_{x,y \in \CoI}\Bigl( \frac{1}{4} \le y \to |x| \le y \to
    \frac{x}{y} \in \CoI \Bigr).
  \end{equation*}
\end{thm}

\begin{proof}
  The proof uses coinduction \eqref{E:GfpI} on $\CoI$ with $P \defeq
  \set{z}{\ex_{x,y \in \CoI}( |x|\le y \land \frac{1}{4} \le y \land
    z=\frac{x}{y})}$.  It suffices to prove \eqref{E:GfpI}'s premise
  for this $P$.  Let $x,y,z$ be given with $x,y \in \CoI$, $|x| \le
  y$, $\frac{1}{4} \le y$ and $z=\frac{x}{y}$.  From $|x| \le y$ we
  have $z \le 1$.  By a trifold application of the closure axiom
  \eqref{E:ClosureI} to $x$ we obtain $d_1, d_2, d_3 \in \mathrm{Sd}$
  and $\tilde{x} \in \CoI$ such that $x =
  \frac{4d_1+2d_2+d_3+\tilde{x}}{8}$ or $x=d_1d_2d_3\tilde{x}$ for
  short.  We now distinguish three cases.

  If $x=1d_2d_3\tilde{x}$, $x=01d_3\tilde{x}$ or $x=001\tilde{x}$,
  then $0 \le x$.  Pick $d=1$ and $z'=\frac{x'}{y}$ with $x' =
  4\frac{x+\frac{-y}{2}}{2}$.  Then $x' \in \CoI$ by
  Lemma~\ref{L:CoIAux}.  From \eqref{E:ShiftIneq} we also obtain $|x'|
  \le y$ and hence $z' \in P$.  One can easily check that
  $z=\frac{1+z'}{2}$.

  If $x=\bar{1}d_2d_3\tilde{x}$, $x=0\bar{1}d_3\tilde{x}$ or
  $x=00\bar{1}\tilde{x}$, then $x \le 0$.  Pick $d=-1$ and
  $z'=\frac{x'}{y}$ with $x' = 4\frac{x+\frac{y}{2}}{2}$.  We can then
  proceed as in the first case.

  The final case is $x=000\tilde{x}$.  Then $|x| \le \frac{1}{8}$;
  pick $d=0$ and $z'=\frac{x'}{y}$ with $x'=2x$.  We obtain $|x'| \le
  \frac{1}{4} \le y$ and with Lemma~\ref{L:Double} also $x' \in \CoI$.
  Therefore $z' \in P$, and $z=\frac{z'}{2}$ is easily checked.
\end{proof}

The term Minlog extracts is displayed in Figure~3.
\begin{figure}[t]
\begin{verbatim}
[u,u0](CoRec ai=>ai)u
 ([u1][case (cCoIClosure u1)
   (s pair u2 -> [case s
     (SdR -> SdR pair InR(cCoIDivSatCoIClAuxR u1 u0))
     (SdM -> [case (cCoIClosure u2)
       (s0 pair u3 -> [case s0
        (SdR -> SdR pair InR(cCoIDivSatCoIClAuxR u1 u0))
        (SdM -> [case (cCoIClosure u3)
         (s1 pair u4 -> [case s1
          (SdR -> SdR pair InR(cCoIDivSatCoIClAuxR u1 u0))
          (SdM -> SdM pair InR(cCoIToCoIDouble u1))
          (SdL -> SdL pair InR(cCoIDivSatCoIClAuxL u1 u0))])])
        (SdL -> SdL pair InR(cCoIDivSatCoIClAuxL u1 u0))])])
     (SdL -> SdL pair InR(cCoIDivSatCoIClAuxL u1 u0))])])
\end{verbatim}
  \caption{Extracted term for Theorem ~\ref{T:CoIDiv}.}
\end{figure}
The three occurrences of \texttt{cCoIClosure} correspond to the
trifold application of the closure axioms \eqref{E:ClosureI} to $x$.  The
seven cases $x=1d_2d_3\tilde{x}$, $x=01d_3\tilde{x}$,
$x=001\tilde{x}$, $x=000\tilde{x}$, $x=\bar{1}d_2d_3\tilde{x}$,
$x=0\bar{1}d_3\tilde{x}$ and $x=00\bar{1}\tilde{x}$ are clearly
visible.  In the first three cases \texttt{SdR} corresponds to picking
$d=1$ and usage of the $\mathrm{AuxR}$ function from
Lemma~\ref{L:CoIAux}, and similarly in the last three cases
\texttt{SdL} corresponds to picking $d=-1$ and usage of
$\mathrm{AuxL}$.  In the middle case \texttt{SdM} corresponds to
$d=0$; there we use the computational content of Lemma~\ref{L:Double}.
We can describe the extracted term by a function $\mathrm{Div} \colon
\typeSds \typeTo \typeSds \typeTo \typeSds$ corecursively defined
by
\begin{equation*}
  \mathrm{Div}(u,v) \defeq
  \begin{cases}
    \mathrm{SdR}::\mathrm{Div}(\mathrm{AuxR}(u,v),v) &\hbox{if
      $u=1\tilde{u} \lor u=01\tilde{u} \lor u=001\tilde{u}$},
    \\ \mathrm{SdM}::\mathrm{Div}(\mathrm{Double}(u),v) &\hbox{if
      $u=000\tilde{u}$},
    \\ \mathrm{SdL}::\mathrm{Div}(\mathrm{AuxL}(u,v),v) &\hbox{if
      $u=\bar{1}\tilde{u} \lor u=0\bar{1}\tilde{u} \lor
      u=00\bar{1}\tilde{u}$}.
  \end{cases}
\end{equation*}
As abbreviation in the case distinction, we have denoted
$\mathrm{SdL},\mathrm{SdM}$ and $\mathrm{SdR}$ by $\overline{1},0$ and
$1$ and we have omitted the constructor $::$.

We use this description of the extracted term to see how far we have
to look into $u$ and $v$ to determine the first $n$ entries of
$\mathrm{Div}(u,v)$.  To this end we write the above equation as
\begin{equation*}
  \mathrm{Div}(u,v) = d(u) :: \mathrm{Div}(G(u,v),v),
\end{equation*}
where $d(u)$ depends on the first three digits of $u$, and $G(u,v)$ is
one of $\mathrm{AuxR}(u,v)$, $\mathrm{Double}(u)$ or
$\mathrm{AuxL}(u,v)$, according to the present case.  Recall
\begin{align*}
  \mathrm{AuxL}(u,v) &\defeq
  \mathrm{Double}(\mathrm{Double}( \mathrm{Av}(u, h(n(v))))),
  \\
  \mathrm{AuxR}(u,v) &\defeq
  \mathrm{Double}(\mathrm{Double}( \mathrm{Av}(u, h(v)))).
\end{align*}
By the equations for $-$, $h$, $\mathrm{Av}$ and $\mathrm{Double}$ we
see that the first $n$ entries of
\begin{alignat*}{2}
  &{-}u &\quad&\hbox{need the first $n$ entries of $u$},
  \\
  &h(u) &\quad&\hbox{need the first $n-1$ entries of $u$},
  \\
  &\mathrm{Av}(u,v) &\quad&
  \hbox{need the first $n+1$ entries of $u$ and $v$},
  \\
  &\mathrm{Double}(u) &\quad&\hbox{need the first $n+1$ entries of $u$}.
\end{alignat*}
Hence $\mathrm{AuxR}(u,v)$, $\mathrm{AuxL}(u,v)$ and $G(u,v)$ all need
at most the first $n+3$ entries of $u$ and $n+2$ entries of $v$.
Iterating the above equation for $G$ gives for $\mathrm{Div}(u,v)$ the
representation
\begin{equation*}
  d(u) :: d(G(u,v)) :: d(G(G(u,v),v)) :: d(G(G(G(u,v),v),v),v) \dots
\end{equation*}
Therefore the first $n$ entries of $\mathrm{Div}(u,v)$ depend on at
most the first $3n$ entries of $u$ and the first $3n-1$ entries of
$v$.

\subsection{Division for binary reflected digit streams}
\label{SS:DivGray}
In the stream representation of real numbers described in
Section~\ref{SS:DivSd} adjacent dyadics of the same length can differ
in many digits:
\begin{equation*}
  \frac{7}{16} \sim 1 \bar{1} 1 1,
  \qquad
  \frac{9}{16} \sim 1 1 \bar{1} \bar{1}.
\end{equation*}
A possible cure is to \emph{flip} after an occurrence of $1$; 
see Figure~4.  The result is called binary reflected (or Gray-) code.
\begin{figure}[ht]
  \begin{tikzpicture}[scale=.16]
    \coordinate [label=below:\textcolor{red}{$0$}] (0) at (32,0);
    \coordinate [label=below:\textcolor{red}{$-\frac{1}{2}$}] (-1/2) at (16,6);
    \coordinate [label=below:\textcolor{red}{$\frac{1}{2}$}] (1/2) at (48,6);
    \coordinate [label=below:\textcolor{red}{$-\frac{3}{4}$}] (-3/4) at (8,12);
    \coordinate (-1/4) at (24,12);
    \coordinate (1/4) at (40,12);
    \coordinate [label=below:\textcolor{red}{$\frac{3}{4}$}] (3/4) at (56,12);
    \coordinate [label=left:\textcolor{red}{$-\frac{7}{8}$}] (-7/8) at (4,18);
    \coordinate (-5/8) at (12,18);
    \coordinate (-3/8) at (20,18);
    \coordinate (-1/8) at (28,18);
    \coordinate (1/8) at (36,18);
    \coordinate (3/8) at (44,18);
    \coordinate (5/8) at (52,18);
    \coordinate [label=right:\textcolor{red}{$\frac{7}{8}$}] (7/8) at (60,18);
    \coordinate[label=left:\textcolor{red}{$-\frac{15}{16}$}](-15/16) at (2,24);
    \coordinate (-13/16) at (6,24);
    \coordinate (-11/16) at (10,24);
    \coordinate (-9/16) at (14,24);
    \coordinate (-7/16) at (18,24);
    \coordinate (-5/16) at (22,24);
    \coordinate (-3/16) at (26,24);
    \coordinate (-1/16) at (30,24);
    \coordinate (1/16) at (34,24);
    \coordinate (3/16) at (38,24);
    \coordinate (5/16) at (42,24);
    \coordinate (7/16) at (46,24);
    \coordinate (9/16) at (50,24);
    \coordinate (11/16) at (54,24);
    \coordinate (13/16) at (58,24);
    \coordinate[label=right:\textcolor{red}{$\frac{15}{16}$}](15/16) at (62,24);
    
    \draw (0) edge node [below] {$\bar{1}$} (-1/2);
    \draw (0) edge node [below] {$1$} (1/2);
    \draw (-1/2) edge node [below] {$\bar{1}$} (-3/4);
    \draw (-1/2) edge node [below] {$1$} (-1/4);

    \draw (1/2) edge node [below] {$1$} (1/4);
    \draw (1/2) edge node [below] {$\bar{1}$} (3/4);
    \draw (-3/4) edge node [left] {$\bar{1}$} (-7/8);
    \draw (-3/4) edge node [right] {$1$} (-5/8);

    \draw (-1/4) edge node [left] {$1$} (-3/8);
    \draw (-1/4) edge node [right] {$\bar{1}$} (-1/8);
    \draw (1/4) edge node [left] {$\bar{1}$} (1/8);
    \draw (1/4) edge node [right] {$1$} (3/8);

    \draw (3/4) edge node [left] {$1$} (5/8);
    \draw (3/4) edge node [right] {$\bar{1}$} (7/8);
    \draw (-7/8) edge node [left] {$\bar{1}$} (-15/16);
    \draw (-7/8) edge node [right] {$1$} (-13/16);

    \draw (-5/8) edge node [left] {$1$} (-11/16);
    \draw (-5/8) edge node [right] {$\bar{1}$} (-9/16);
    \draw (-3/8) edge node [left] {$\bar{1}$} (-7/16);
    \draw (-3/8) edge node [right] {$1$} (-5/16);

    \draw (-1/8) edge node [left] {$1$} (-3/16);
    \draw (-1/8) edge node [right] {$\bar{1}$} (-1/16);
    \draw (1/8) edge node [left] {$\bar{1}$} (1/16);
    \draw (1/8) edge node [right] {$1$} (3/16);

    \draw (3/8) edge node [left] {$1$} (5/16);
    \draw (3/8) edge node [right] {$\bar{1}$} (7/16);
    \draw (5/8) edge node [left] {$\bar{1}$} (9/16);
    \draw (5/8) edge node [right] {$1$} (11/16);

    \draw (7/8) edge node [left] {$1$} (13/16);
    \draw (7/8) edge node [right] {$\bar{1}$} (15/16);
  \end{tikzpicture}
  \caption{Gray code.}
\end{figure}

Then two dyadics of the same length are adjacent if and only if they differ
in exactly one digit.  For instance we now have
\begin{equation*}
  \frac{7}{16} \sim 1 1  1 \bar{1},
  \qquad
  \frac{9}{16} \sim 1 \bar{1} 1 \bar{1}.
\end{equation*}
This is a desirable property of stream-coded real numbers.

Again when doing arithmetical operations on Gray code the problem of
productivity arises, which is dealt with in a somewhat different way.
The idea
is to introduce two \inquotes{modes} when generating the code, and
flip from one mode to the other whenever we encounter the digit $1$.
More precisely, instead of the predicate $I$ we now use two predicates
$G,H$ ($G$ for Gray, $H$ next character) and flip from one to the other
after reading $1$.  They are defined by simultaneous induction, with
clauses
\begin{alignat*}{2}
  &\forall_{d \in \mathrm{Psd}} \forall_{x \in G}\Bigl(
  -d\frac{x-1}{2} \in G \Bigr),
  &\quad&
  \forall_{x \in H}\Bigl( \frac{x}{2} \in G \Bigr),
  \\
  &\forall_{d \in \mathrm{Psd}} \forall_{x \in G}\Bigl(
  d\frac{x+1}{2} \in H \Bigr),
  &\quad&
  \forall_{x \in H}\Bigl( \frac{x}{2} \in H \Bigr).
\end{alignat*}
Here $\mathrm{Psd}$ (for proper signed digit) is the (inductive)
predicate consisting of the integers $1,-1$.  We require that
$\mathrm{Psd}$ has computational content, in the booleans $\typeB =
\{\true, \false\}$.  Note that there are no nullary clauses (as for
$I$), since they are not needed.  We are only interested in the duals
$\CoG, \CoH$, coinductively defined by the simultaneous closure axioms
\begin{equation}
  \label{E:ClosureGH}
  \begin{split}
    &\forall_{x \in \CoG}\Bigl( 
    \ex_{d \in \mathrm{Psd}} \ex_{x' \in \CoG} \Bigl( x=-d\frac{x'-1}{2} \Bigr) \lor
    \ex_{x' \in \CoH} \Bigl( x=\frac{x'}{2} \Bigr) \Bigr)
    \\
    &\forall_{x \in \CoH}\Bigl( 
    \ex_{d \in \mathrm{Psd}} \ex_{x' \in \CoG} \Bigl( x=d\frac{x'+1}{2} \Bigr) \lor
    \ex_{x' \in \CoH} \Bigl( x=\frac{x'}{2} \Bigr) \Bigr)
  \end{split}
\end{equation}
and the simultaneous coinduction axioms
\begin{equation}
  \label{E:GfpGH}
  \begin{split}
    &\forall_{x \in P}\Bigl( 
    \ex_{d \in \mathrm{Psd}}\ex_{x' \in \CoG \cup P} \Bigl( x=-d\frac{x'-1}{2} \Bigr)
    \lor
    \ex_{x' \in \CoH \cup Q} \Bigl( x=\frac{x'}{2} \Bigr) \Bigr) \to {}
    \\
    &\forall_{x \in Q}\Bigl( 
    \ex_{d \in \mathrm{Psd}} \ex_{x' \in \CoG \cup P} \Bigl( x=d\frac{x'+1}{2} \Bigr)
    \lor
    \ex_{x' \in \CoH \cup Q} \Bigl( x=\frac{x'}{2} \Bigr) \Bigr) \to {}
    \\[3pt]
    &P \subseteq \CoG
    \qquad \hbox{and the same with conclusion $Q \subseteq \CoH$}.
  \end{split}
\end{equation}

We require that the predicates $G,H$ and therefore $\CoG, \CoH$ as well
have computational content, in simultaneously defined algebras
$\typeG$ and $\typeH$ given by the constructors
\begin{alignat*}{3}
  &\constrLr \colon \typeB \typeTo \typeG \typeTo \typeG,
  &\qquad&
  \constrU \colon \typeH \typeTo \typeG
  &\qquad&
  \hbox{for $\typeG$},
  \\
  &\constrFin \colon \typeB \typeTo \typeG \to \typeH,
  &\qquad&
  \constrD \colon \typeH \typeTo \typeH
  &\qquad&
  \hbox{for $\typeH$}.
\end{alignat*}
We write $\constrLr_\true(u)$ for $\constrLr(\true, u)$ and
$\constrLr_\false(u)$ for $\constrLr(\false, u)$, and similarly for
$\constrFin$.  Then $\constrLr_{\true/\false}$ means left/right,
$\constrFin_{\true/\false}$ means finally left/right.  The delay
constructors are $\constrU$ (\inquotes{undefined}) for $\typeG$ and
$\constrD$ (\inquotes{delay}) for $\typeH$.  Figure~5 indicates
how these constructors generate Gray code.
\begin{figure}[ht]
  \begin{tikzpicture}[scale=.22]
    \coordinate [label=below:\textcolor{red}{$0$}] (0) at (32,0);
    \coordinate [label=below:\textcolor{red}{$-\frac{1}{2}$}] (-1/2) at (16,6);
    \coordinate [label=below:\textcolor{red}{$\frac{1}{2}$}] (1/2) at (48,6);
    \coordinate [label=below:\textcolor{red}{$-\frac{3}{4}$}] (-3/4) at (8,12);
    \coordinate (-1/4) at (24,12);
    \coordinate (1/4) at (40,12);
    \coordinate (3/4) at (56,12);
    \coordinate (-7/8) at (4,18);
    \coordinate [label=below:\textcolor{red}{$\frac{3}{4}$}] (3/4) at (56,12);
    \coordinate (-5/8) at (12,18);
    \coordinate (-3/8) at (20,18);
    \coordinate (-1/8) at (28,18);
    \coordinate (1/8) at (36,18);
    \coordinate (3/8) at (44,18);
    \coordinate (5/8) at (52,18);
    \coordinate (7/8) at (60,18);

    \draw [fill] (32,0) circle [radius=0.5];
    \draw [fill] (16,6) circle [radius=0.5];
    \draw [thick] (32,6) circle [radius=0.5];
    \draw [fill] (48,6) circle [radius=0.5];
    \draw [fill] (8,12) circle [radius=0.5];
    \draw [thick] (16,12) circle [radius=0.5];
    \draw [fill] (24,12) circle [radius=0.5];
    \draw [thick] (32,12) circle [radius=0.5];
    \draw [fill] (40,12) circle [radius=0.5];
    \draw [thick] (48,12) circle [radius=0.5];
    \draw [fill] (56,12) circle [radius=0.5];
    \draw [fill] (4,18) circle [radius=0.5];
    \draw [thick] (8,18) circle [radius=0.5];
    \draw [fill] (12,18) circle [radius=0.5];
    \draw [thick] (16,18) circle [radius=0.5];
    \draw [fill] (20,18) circle [radius=0.5];
    \draw [thick] (24,18) circle [radius=0.5];
    \draw [fill] (28,18) circle [radius=0.5];
    \draw [thick] (32,18) circle [radius=0.5];
    \draw [fill] (36,18) circle [radius=0.5];
    \draw [thick] (40,18) circle [radius=0.5];
    \draw [fill] (44,18) circle [radius=0.5];
    \draw [thick] (48,18) circle [radius=0.5];
    \draw [fill] (52,18) circle [radius=0.5];
    \draw [thick] (56,18) circle [radius=0.5];
    \draw [fill] (60,18) circle [radius=0.5];

    \draw (0) edge node [below] {$\constrLr_\false$} (-1/2);
    \draw (0) edge node [right] {$\constrU$} (32,6);
    \draw (0) edge node [below] {$\constrLr_\true$} (1/2);
    \draw (-1/2) edge node [below] {$\constrLr_\false$} (-3/4);
    \draw (-1/2) edge node [right] {$\constrU$} (16,12);
    \draw (-1/2) edge node [right] {$\constrLr_\true$} (24,12);
    \draw [dashed] (32,6) edge node [left] {$\constrFin_\false$} (24,12);
    \draw [dashed] (32,6) edge node [right] {$\constrD$} (32,12);
    \draw [dashed] (32,6) edge node [right] {$\constrFin_\true$} (40,12);
    \draw (1/2) edge node [below] {$\constrLr_\true$} (40,12);
    \draw (1/2) edge node [right] {$\constrU$} (48,12);
    \draw (1/2) edge node [below] {$\constrLr_\false$} (3/4);
    \draw (-3/4) edge (4,18);
    \draw (-3/4) edge (8,18);
    \draw (-3/4) edge (12,18);
    \draw [dashed] (16,12) edge (12,18);
    \draw [dashed] (16,12) edge (16,18);
    \draw [dashed] (16,12) edge (20,18);
    \draw (24,12) edge (20,18);
    \draw (24,12) edge (24,18);
    \draw (24,12) edge (28,18);
    \draw [dashed] (32,12) edge (28,18);
    \draw [dashed] (32,12) edge (32,18);
    \draw [dashed] (32,12) edge (36,18);
    \draw (40,12) edge (36,18);
    \draw (40,12) edge (40,18);
    \draw (40,12) edge (44,18);
    \draw [dashed] (48,12) edge (44,18);
    \draw [dashed] (48,12) edge (48,18);
    \draw [dashed] (48,12) edge (52,18);
    \draw (56,12) edge (52,18);
    \draw (56,12) edge (56,18);
    \draw (56,12) edge (60,18);
  \end{tikzpicture}
  \caption{Gray code with delay.}
\end{figure}

Totality $\Total_\typeG$, $\Total_\typeH$ and cototality
$\CoTotal_\typeG$, $\CoTotal_\typeH$ can be defined as for $\typeIntv$
above.  Objects $u \in \CoTotal_\typeG$ are called \emph{streams} in
Gray code.

The computational content of the closure axioms \eqref{E:ClosureGH}
are the destructors $\destr_\typeG \colon \typeG \typeTo (\typeB
\typeProd \typeG) \typeSum \typeH$ and $\destr_\typeH \colon \typeH
\typeTo (\typeB \typeProd \typeG) \typeSum \typeH$ defined by
\begin{alignat*}{2}
  &\destr_\typeG(\constrLr_b(u)) = \termSumIntroLeft \pair {b} {u},
  &\quad&
  \destr_\typeG(\constrU v) = \termSumIntroRight(v),
  \\
  &\destr_\typeH(\constrFin_b(u)) = \termSumIntroLeft \pair {b} {u},
  &\quad&
  \destr_\typeH(\constrD v) = \termSumIntroRight(v).
\end{alignat*}
The computational content of the coinduction axioms \eqref{E:GfpGH}
are instances of the simultaneous corecursion operators
\begin{align}
  \label{E:SimCoRecType}
  \begin{split}
    &\corec_\typeG^{(\typeG, \typeH), (\sigma,\tau)} \colon \sigma
    \typeTo \delta_\typeG \typeTo \delta_\typeH \typeTo \typeG
    \\
    &\corec_\typeH^{(\typeG, \typeH), (\sigma,\tau)} \colon \tau
    \typeTo \delta_\typeG \typeTo \delta_\typeH \typeTo \typeH
  \end{split}
\end{align}
with step types
\begin{align*}
  &\delta_\typeG \defeq \sigma \typeTo \typeB \typeProd( \typeG
  \typeSum \sigma) \typeSum (\typeH \typeSum \tau),
  \\ &\delta_\typeH \defeq \tau \typeTo \typeB \typeProd( \typeG
  \typeSum \sigma) \typeSum (\typeH \typeSum \tau).
\end{align*}
The type $\typeB \typeProd (\typeG \typeSum \sigma) \typeSum (\typeH
\typeSum \tau)$ appears since $\typeG$ has the two constructors
$\constrLr \colon \typeB \typeTo \typeG \typeTo \typeG$ and $\constrU
\colon \typeH \typeTo \typeG$, and $\typeH$ has the two constructors
$\constrFin \colon \typeB \typeTo \typeG \typeTo \typeH$ and $\constrD
\colon \typeH \typeTo \typeH$.  Omitting the upper indices of $\corec$
the computation rules for the terms $\corec_\typeG s g h$ and
$\corec_\typeH t g h$ are
\begin{equation}
  \label{E:GrayDivShiftGH}
  \begin{split}
    \corec_\typeG s g h =
    \begin{cases}
      \constrLr(b,u)
      &\hbox{if $\eqd{g(s)}{\termSumIntroLeft\pair{b}
          {\termSumIntroLeft_{\typeG \typeSum \sigma}u}}$}
      \\
      \constrLr(b,\corec_\typeG s' g h)
      &\hbox{if $\eqd{g(s)}{\termSumIntroLeft\pair{b}
          {\termSumIntroRight_{\typeG \typeSum \sigma}s'}}$}
      \\
      \constrU(v)
      &\hbox{if $\eqd{g(s)}{\termSumIntroRight
          (\termSumIntroLeft_{\typeH \typeSum \tau}v})$}
      \\
      \constrU(\corec_\typeH t g h)
      &\hbox{if $\eqd{g(s)}{\termSumIntroRight
          (\termSumIntroRight_{\typeH \typeSum \tau}t})$}
    \end{cases}
    \\
    \corec_\typeH t g h =
    \begin{cases}
      \constrFin(b,u)
      &\hbox{if $\eqd{h(t)}{\termSumIntroLeft\pair{b}
          {\termSumIntroLeft_{\typeG \typeSum \sigma}u}}$}
      \\
      \constrFin(b,\corec_\typeG s g h)
      &\hbox{if $\eqd{h(t)}{\termSumIntroLeft\pair{b}
          {\termSumIntroRight_{\typeG \typeSum \sigma}s}}$}
      \\
      \constrD(v)
      &\hbox{if $\eqd{h(t)}{\termSumIntroRight
          (\termSumIntroLeft_{\typeH \typeSum \tau}v})$}
      \\
      \constrD(\corec_\typeH t' g h)
      &\hbox{if $\eqd{h(t)}{\termSumIntroRight
          (\termSumIntroRight_{\typeH \typeSum \tau}t'})$}
    \end{cases}
  \end{split}
\end{equation}
with $s$ of type $\sigma$ and $t$ of type $\tau$.

Our goal in this section is to prove that $[-1,1]$ is closed under
division $\frac{x}{y}$ (for $|x| \le |y| > 0$), w.r.t.\ the
representation of reals in Gray code.  We proceed essentially as for
the signed digit case.  The main difference is that the simultaneous
definition of $\CoG$ and $\CoH$ makes it necessary to use simultaneous
coinduction.

\begin{thm}
  \label{T:CoGAverage}
  The average of two real numbers $x,y$ in $\CoG$ is in $\CoG$:
  \begin{equation*}
    \forall_{x,y \in \CoG}\Bigl( \frac{x+y}{2} \in \CoG \Bigr).
  \end{equation*}
\end{thm}

\begin{proof}
  The proof in \cite{BergerMiyamotoSchwichtenbergTsuiki16} can be
  adapted to the present setting with concrete rather than abstract
  reals.
\end{proof}

\begin{lem}
  \label{L:CoGCoHMinus}
  $\CoG$ and $\CoH$ are closed under minus:
  \begin{equation*}
    \forall_{x\in \CoG}(-x\in \CoG),
    \qquad
    \forall_{x\in \CoH}(-x\in \CoH).
  \end{equation*}
\end{lem}

\begin{proof}
  We prove both claims simultaneously, by coinduction:
  \begin{align*}
    \begin{split}
      &\forall_{x \in P}\Bigl( 
      \ex_{d \in \mathrm{Psd}} \ex_{x' \in \CoG \cup P} \Bigl(
      x=-d\frac{x'-1}{2} \Bigr)
      \lor
      \ex_{x' \in \CoH \cup Q} \Bigl( x=\frac{x'}{2} \Bigr) \Bigr) \to {}
      \\
      &\forall_{x \in Q}\Bigl( 
      \ex_{d \in \mathrm{Psd}} \ex_{x' \in \CoG \cup P} \Bigl(
      x=d\frac{x'+1}{2} \Bigr)
      \lor
      \ex_{x' \in \CoH \cup Q} \Bigl( x=\frac{x'}{2} \Bigr) \Bigr) \to {}
      \\[3pt]
      &\forall_x(x\in P \to x \in \CoG) \land
      \forall_x(x\in Q \to x \in \CoH).
    \end{split}
  \end{align*}
  We choose $P \defeq \set{x}{-x \in \CoG}$ and $Q \defeq \set{x}{-x
    \in \CoH}$. Since $\CoG$ and $\CoH$ are invariant under real
  equality (also proven by coinduction) we have to prove the two
  premises.  We only consider the first one; the second is proved
  similarly.  Let $x \in P$ be given.  The goal is
  \begin{equation*}
    \ex_{d \in \mathrm{Psd}} \ex_{x' \in \CoG \cup P} \Bigl(
    x=-d\frac{x'-1}{2} \Bigr) \lor
    \ex_{x' \in \CoH \cup Q} \Bigl( x=\frac{x'}{2} \Bigr).
  \end{equation*}
  From $x \in P$ we get $-x\in \CoG$ and therefore 
  \begin{equation*}
    \ex_{d \in \mathrm{Psd}} \ex_{x' \in \CoG} \Bigl(
    -x=-d\frac{x'-1}{2} \Bigr) \lor
    \ex_{x' \in \CoH} \Bigl( -x=\frac{x'}{2} \Bigr)
  \end{equation*}
  by the closure axiom.  If the right hand side of this disjunction
  holds, we are done by definition of $Q$.  If the left hand side
  holds, we have $x=d\frac{x'-1}{2}$ for some $d\in \mathrm{Psd}$ and
  $x' \in \CoG$.  To get the left hand side of the disjunction in the
  goal formula, we take $d:=-d$ and $x':=x'$.
\end{proof}

The function we get from the proof takes a Gray code and flips all
$\true$ to $\false$ and the other way around.  We denote this function
by $-$.  Then the computation rules are:
\begin{equation*}
  \begin{split}
    -(\constrLr_\true(u)) &= \constrLr_\false(-u),
    \\
    -(\constrLr_\false(u)) &= \constrLr_\true( -u),
    \\
    -(\constrU(v)) &= \constrU(-v),
  \end{split}
  \quad
  \begin{split}
    -(\constrFin_\true(u)) &= \constrFin_\false(-u),
    \\
    -(\constrFin_\false(u)) &= \constrFin_\true( -u)),
    \\
    -(\constrD(v)) &= \constrD(-v) .
  \end{split}
\end{equation*}

\begin{lem}
  \label{L:CoGEquiCoH}
  $\CoG$ and $\CoH$ are equivalent:
  \begin{equation*}
    \forall_{x\in \CoG}(x\in \CoH),
    \\\qquad
    \forall_{x\in \CoH}(x\in \CoG).
  \end{equation*}
\end{lem}

\begin{proof}
  In the coinduction axiom we choose $P \defeq \CoH$ and $Q \defeq
  \CoG$ and the conclusion in the coinduction axiom is our goal.
  Therefore we have to prove the premises of the coinduction axiom.
  We only consider the first one; the secound one is proved similarly.
  Let $x \in \CoH$ be given.  By the closure axiom of $\CoH$ we have
  \begin{equation*}
    \ex_{d \in \mathrm{Psd}} \ex_{x' \in \CoG} \Bigl(
    x=d\frac{x'+1}{2} \Bigr) \lor
    \ex_{x' \in \CoH} \Bigl( x=\frac{x'}{2} \Bigr)
  \end{equation*}
  If the second part of this disjunction holds, we get directly the
  second part of the first premise.  If the first part of the
  disjunction holds, we have $x=d\frac{x'+1}{2}$ for some $d \in
  \mathrm{Psd}$ and $x'\in \CoG$.  Then $x=-d\frac{(-x')-1}2$ and
  using Lemma~\ref{L:CoGCoHMinus} we get the first part of disjunction
  in the premise.
\end{proof}

As computational content we have two functions
$\mathrm{ToCoH}\colon\mathbb{G}\to \mathbb{H}$ and
$\mathrm{ToCoG}\colon\mathbb{H}\to\mathbb{G}$, with computation rules
\begin{equation*}
  \begin{split}
    \mathrm{ToCoH}(\constrLr_b(u)) &= \constrFin_b(-u),
    \\
    \mathrm{ToCoH}(\constrU(v)) &= \constrD(v),
  \end{split}
  \quad
  \begin{split}
    \mathrm{ToCoG}(\constrFin_b(u)) &= \constrLr_b(-u),
    \\
    \mathrm{ToCoG}(\constrD(v)) &= \constrU(v) .
  \end{split}
\end{equation*}
Note that no corecursive call is involved.  We denote
$\mathrm{ToCoG}(v)$, $\mathrm{ToCoH}(u)$ by $\tilde{v}$, $\tilde{u}$,
respectively.

We show that $\CoG$ is closed under shifting a real $x\leq 0$ ($x \geq
0$) by $+1$ ($-1$):
  \begin{equation*}
    \forall_{x \in \CoG}(( x \le 0 \to x+1 \in \CoG)\land
    ( 0 \le x \to x-1 \in \CoG)).
  \end{equation*}
The proof is by coinduction, simultaneously with the same proposition
for $\CoH$.  For the coinduction proof it is helpful to reformulate the
claim in such a way that the conclusion has the form $x \in \CoG$:

\begin{lem}
  \label{L:CoGShift}
  \begin{equation*}
    \forall_x(\ex_{y\in \CoG}( y\le 0 \land (x=y+1 \lor x=-(y+1))) \to
    x  \in \CoG).
  \end{equation*}
\end{lem}

\begin{proof}
  We use the introduction axiom for $\CoG$ and $\CoH$, and
  simultaneously prove
  \begin{equation*}
    \forall_x(\ex_{y\in \CoH}(y\leq 0 \land (x=y+1 \lor x=-(y+1))) \to
    x \in \CoH).
  \end{equation*}
  The competitor predicates are 
  \begin{align*}
    P &\defeq \set{x}{\ex_{y \in \CoG}( y \le 0 \land (x=y+1 \lor x=-(y+1)))},
    \\
    Q &\defeq \set{x}{\ex_{y \in \CoH}( y \le 0 \land (x=y+1 \lor x=-(y+1)))}.
  \end{align*}
  It suffices to prove the premises of the coinduction axiom for these
  $P$, $Q$.  We only consider the first part; the second is proved
  similarly.  Let $x\in P$ be given.  Then we have $y$ with $y \in
  \CoG$, $y \le 0$ and $x=y+1 \lor x=-(y+1)$.  From $y \in \CoG$ we
  know $|y| \le 1$ and with $y \le 0$ also $|x| \le 1$.  We need to
  prove the disjunction
  \begin{equation*}
    D \defeq \ex_{d \in \mathrm{Psd}} \ex_{x' \in \CoG \cup P}\Bigl(
    x = -d\frac{x'-1}{2} \Bigr) \lor
    \ex_{x' \in \CoH \cup Q}\Bigl(x = \frac{x'}{2} \Bigr).
  \end{equation*}
  From $y \in \CoG$ we obtain either (i) $e \in \mathrm{Psd}$, $z \in
  \CoG$ with $y = -e\frac{z-1}{2}$ or else (ii) $z \in \CoH$ with $y =
  \frac{z}{2}$.  In case (i) we have $|z| \le 1$ since $z \in \CoG$.
  We distinguish cases on $e \in \mathrm{Psd}$, and use
  Lemma~\ref{L:CoGCoHMinus}.

  \emph{Case} $e=1$.  Then $0 \ge y = -\frac{z-1}{2} \ge -\frac{1}{2}
  + \frac{1}{2} = 0$, hence $y=0$.  If $x=y+1$, take $x=1$ and go for
  the l.h.s.\ of the disjunction $D$.  Picking $d=1$ and $x'= -1$
  gives the claim (since $-1 \in \CoG$).  If $x= -(y+1)$, take $x= -1$
  and go for the l.h.s.\ of the disjunction $D$.  Picking $d= -1$ and
  $x'= -1$ gives the claim (since $-1 \in \CoG$).

  \emph{Case} $e= -1$.  If $x=y+1$, go for the l.h.s.\ of the
  disjunction $D$, picking $d=1$ and $x'= -z$.  Then $x = y+1 =
  \frac{z-1}{2} +1 = \frac{z+1}{2} = -\frac{-z-1}{2}$.  If $x=
  -(y+1)$, go for the l.h.s.\ of the disjunction $D$, picking $d= -1$
  and $x'= -z$.  Then $x = -(y+1) = -\frac{z-1}{2} -1 = \frac{-z+1}{2}
  -1 = \frac{-z-1}{2}$.

  In case (ii) we have $|z| \le 1$ (since $z \in \CoH$) and $z \le 0$
  (since $y = \frac{z}{2}$ and $y \le 0$).  Hence $|z+1| \le 1$.  If
  $x=y+1$, go for the l.h.s.\ of the disjunction $D$, picking $d=1$
  and $x'= -(z+1)$.  Then $-(z+1) \in P$ (since also $z \in \CoG$ and
  $z \le 0$), and $x = y+1 = \frac{z}{2} +1 = \frac{z+2}{2} =
  -\frac{-(z+1)-1}{2}$.  If $x=-(y+1)$, go for the l.h.s.\ of the
  disjunction $D$, with $d= -1$ and $x'= -(z+1)$.  Then again $-(z+1)
  \in P$ and $x = -(y+1) = \frac{-(z+1)-1}{2}$.
\end{proof}

As computational content of a formalization of this proof Minlog
returns a term involving simultaneous corecursion.  The
corresponding algorithm can be described as a pair of two functions $\shG
\colon \typeG \typeTo \typeB \typeTo \typeG$ and $\shH \colon \typeH
\typeTo \typeB \typeTo \typeH$ defined by
\begin{equation*}
  \begin{split}
    \shG(\constrLr_\true(u), b) &= \constrLr_b( -1),
    \\
    \shG(\constrLr_\false(u), b) &= \constrLr_b( -u),
    \\
    \shG(\constrU(v), b) &= \constrLr_b( \shH(\tilde{v}, \false)),
  \end{split}
  \quad
  \begin{split}
    \shH(\constrFin_\true(u), b) &= \constrFin_b(1),
    \\
    \shH(\constrFin_\false(u), b) &= \constrFin_b( -u)),
    \\
    \shH(\constrD(v), b) &= \constrFin_b( \shG(\tilde{v}, \true)).
  \end{split}
\end{equation*}
Recall that $v \mapsto \tilde{v}$ is the function extracted from
the proof of Lemma~\ref{L:CoGEquiCoH} stating that $\CoG$ and $\CoH$
are equivalent.  If $b=\true$, we have computational content of the
first part of the lemma.  If $b=\false$, we have the second part.

\begin{lem}
  \label{L:CoGDouble}
  For $x$ in $\CoG$ with $|x| \le \frac{1}{2}$ we have $2x$ in $\CoG$:
  \begin{equation*}
    \forall_{x \in \CoG}\Bigl( |x| \le \frac{1}{2} \to 2x \in \CoG \Bigr).
  \end{equation*}
\end{lem}

\begin{proof}
  Let $x \in \CoG$ be given.  From the closure axiom for $\CoG$ we
  obtain either (i) $d \in \mathrm{Psd}$, $x' \in \CoG$ with $x=
  -d\frac{x'-1}{2}$, or else (ii) $x' \in \CoH$ with $x=\frac{x'}{2}$.
  In case (i) we distinguish cases on $d \in \mathrm{Psd}$.  In case
  $d=1$ we obtain $x=\frac{-x'+1}{2}$ and hence $2x=-x'+1$.  Since
  $|x| \le \frac{1}{2}$ we have $x' \le 0$.  Now the first part of
  Lemma~\ref{L:CoGShift} gives the claim.  The case $d=-1$ is similar,
  using the second part of Lemma~\ref{L:CoGShift}.  In case (ii) we
  have $2x=x'$, hence the goal $2x \in \CoG$ follows from the
  equivalence of $\CoG$ and $\CoH$.
\end{proof}

Using arguments similar to those in the remark after
Lemma~\ref{L:CoGShift} we can see that the corresponding algorithm can
be written as a function $\mathrm{Double} \colon \typeG \typeTo
\typeG$ with
\begin{align*}
  \mathrm{Double}(\constrLr_\true(u)) &\defeq \shG(-u, \true),
  \\
  \mathrm{Double}(\constrLr_\false(u)) &\defeq \shG(-u, \false),
  \\
  \mathrm{Double}(\constrU(v)) &\defeq \tilde{v}
\end{align*}
where $\shG$ and $v \mapsto \tilde{v}$ are the functions from this
remark.

Parallel to Lemma~\ref{L:CoIAux} we can now prove


\begin{lem}
  \label{L:CoGAux}
  For $x,y$ in $\CoG$ with $\frac{1}{4} \le y$, $|x| \le y$ and
  $0 \le x$ ($x \le 0$) we have $2x-y$ (\,$2x+y$) in $\CoG$:
  \begin{align*}
    &\forall_{x,y \in \CoG} \Bigl( \frac{1}{4} \le y \to |x| \le y \to 0\le x \to
    4\frac{x+\frac{-y}{2}}{2} \in \CoG \Bigr),
    \\
    &\forall_{x,y \in \CoG} \Bigl(\frac{1}{4} \le y \to |x| \le y \to x\le 0 \to
    4\frac{x+\frac{y}{2}}{2} \in \CoG \Bigr).
  \end{align*}
\end{lem}

\begin{proof}
  We only consider the first claim, and again use
  Theorem~\ref{T:CoGAverage} on the average.  In the formulas above
  instead of $2x \pm y$ we have written $4\frac{x+(\pm y/2)}{2}$ to
  make Theorem~\ref{T:CoGAverage} applicable.  One can see easily that
  $\CoG$ is closed under $x \mapsto \frac{x}{2}$.  Formally, this is
  proved by coinduction, and the corresponding algorithm is given by a
  function $h \colon \typeG \typeTo \typeG$.  To prove the first
  claim, let $x,y$ in $\CoG$ with $\frac{1}{4} \le y$, $|x| \le y$ and
  $0 \le x$.  Then clearly $\frac{-y}{2} \in \CoG$, and
  $\frac{x+\frac{-y}{2}}{2} \in \CoG$ by Theorem~\ref{T:CoGAverage}.
  We have the same estimate \eqref{E:ShiftIneq} as in the proof of
  Lemma~\ref{L:CoIAux}; hence we can apply Lemma~\ref{L:CoGDouble}
  twice and obtain $4\frac{x+\frac{-y}{2}}{2} \in \CoG$.  The proof of
  the second claim is similar.  
\end{proof}

The computational content are functions $\mathrm{AuxL}, \mathrm{AuxR}
\colon \typeG \typeTo \typeG \typeTo \typeG$ with
\begin{align*}
  \mathrm{AuxL}(u,u') &\defeq
  \mathrm{Double}(\mathrm{Double}( \mathrm{Av}(u, h(-u')))),
  \\
  \mathrm{AuxR}(u,u') &\defeq
  \mathrm{Double}(\mathrm{Double}( \mathrm{Av}(u, h(u')))).
\end{align*}

Parallel to the essential part of Theorem~\ref{T:CoIDiv} we can now
prove that division by $y$ satisfies to closure axiom for $\CoI$:

\begin{thm}
  \label{T:CoGDivSatCoICl}
  For $x,y$ in $\CoG$ with $\frac{1}{4} \le y$ and $|x| \le y$ we find
  a signed digit $d$ such that $2x=x'+yd$ for some  $x' \in \CoG$
  with $|x'| \le y$:
  \begin{equation*}
    \forall_{x,y \in \CoG} \Bigl( \frac{1}{4} \le y \to |x| \le y \to
    \ex_{d \in \mathrm{Sd}} \ex_{x' \in \CoG} \Bigl( |x'| \le y \land
    \frac{x}{y} = \frac{\frac{x'}{y} + d}{2}\Bigr) \Bigr).
  \end{equation*}
\end{thm}

\begin{proof}
  We proceed as for Theorem~\ref{T:CoIDiv}, now using
  Lemma~\ref{L:CoGAux}.  The proof uses a trifold application of the
  coinduction axioms for $\CoG$ and $\CoH$ to obtain the first three
  digits.
\end{proof}

The computational content is a function $f \colon \typeG
\typeTo \typeG \typeTo \typeSd \typeProd \typeG$ defined by
\begin{align*}
  f(\constrLr_\true(u), u') &= (1, \mathrm{AuxR}(\constrLr_\true(u), u'))
  \\
  f(\constrLr_\false(u), u') &= (-1, \mathrm{AuxL}(\constrLr_\false(u), u'))
  \\
  f(\constrU(v), \constrFin_\true(u)) &=
  (1, \mathrm{AuxR}(\constrU(v), \constrFin_\true(u)))
  \\
  f(\constrU(v), \constrFin_\false(u)) &=
  (-1, \mathrm{AuxL}(\constrU(v), \constrFin_\false(u)))
  \\
  f(\constrU(v), \constrD(\constrFin_\true(u))) &=
  (1, \mathrm{AuxR}(\constrU(v), \constrD(\constrFin_\true(u))))
  \\
  f(\constrU(v), \constrD(\constrFin_\false(u))) &=
  (-1, \mathrm{AuxL}(\constrU(v), \constrD(\constrFin_\false(u))))
  \\
  f(\constrU(v), \constrD(\constrD(v'))) &= (0, \constrD(\constrU(v)))
\end{align*}

\begin{cor}
  \label{C:CoGDiv}
  For $x,y$ in $\CoG$ with $\frac{1}{4} \le y$ and $|x| \le y$ we have
  $\frac{x}{y}$ in $\CoG$:
  \begin{equation*}
    \forall_{x,y \in \CoG}\Bigl( \frac{1}{4} \le y \to |x| \le y \to
    \frac{x}{y} \in \CoG \Bigr).
  \end{equation*}
\end{cor}

\begin{proof}
  By coinduction, simultaneously with the same formula where $\CoG$ is
  replaced by $\CoH$.  Theorem~\ref{T:CoGDivSatCoICl} is used in both
  cases of the simultaneous coinduction.  The extracted term could be
  analyzed in a similar way as in the remark after
  Theorem~\ref{T:CoIDiv}.
\end{proof}

\subsection{Translation to Haskell}
\label{SS:Haskell}
The terms extracted from Theorem~\ref{T:CoIDiv} and
Corollary~\ref{C:CoGDiv} can be translated into Scheme or Haskell
programs.  Because of the presence of
corecursion operators in the extracted terms the lazy evaluation of
Haskell is more appropriate.  As tests we have run (in \texttt{ghci}
with time measuring by \texttt{:set +s}) the signed digit division
files\footnote{\texttt{sddiv.scm} and \texttt{graydiv.scm} residing in
  \texttt{minlog/examples/analysis}} on approximations of
$\frac{1}{3}$ and $\frac{1}{2}$.  To return the first 19 digits of the
result of dividing $\frac{1001}{3001}$ by $\frac{10001}{20001}$ took
about 0.04 seconds in the signed digit case and about 0.06 seconds in
the Gray case.

In the following table we see how the runtime increases if we increase
the number of signed digits in the output.  Of course, this depends on
the used computer, but instead of the concrete numbers we are
interested in the scale of the runtime.

\begin{figure}[ht]
  \begin{center}
    \begin{tabular}{d{5}d{4.2}}
      \toprule
      \multicolumn{1}{c}{\bfseries number of digits} & \multicolumn{1}{c}{\bfseries runtime in seconds} \\
      \midrule
      10 & 0.01\\
      25 & 0.05 \\
      50 & 0.14\\
      75 & 0.26\\
      100& 0.46\\
      250 & 2.69\\
      500 & 10.11\\
      750 & 23.90\\
      1000 & 42.50\\
      2000& 182.92\\
      10000& 4567.74\\
      \bottomrule
    \end{tabular}
  \end{center}
  \caption{Runtime}
\end{figure}

In the remark after Theorem \ref{T:CoIDiv} we have shown that the
look-ahead of the input is linear in the digits of the output, i.e., we
need at most the first $3n$ entries of $u$ and $v$ to compute the
first $n$ entries of $\text{Div}(u,v)$.  But here we see that the
runtime is clearly not linear. This is not surprising because in this
remark we had the representation
\begin{equation*}
  d(u) :: d(G(u,v)) :: d(G(G(u,v),v)) :: d(G(G(G(u,v),v),v),v) \dots
\end{equation*}
of $\text{Div}(u,v)$.  Therefore, for the first $n$ digits, we have to
compute
\begin{equation*}
  G(u,v),G(G(u,v),v),\dots,\underbrace{G(\dots G}_{n-1}(u,v)\dots)
\end{equation*}
and we have to read the first $3n$ digits of $u$ and $v$ and operate
on them $n-1$ times.  We see that $n$ occurs twice in the calculation
and hence the runtime has to be at least quadratic in the numbers of
digits of the output.  We also see this in the table above:

If we compare, for example, the runtimes for 100, 1000 and 10000
digits, we see that a multiplication of the numbers of digits by 10
causes a multiplication of the runtime by approximately 100.
Therefore the runtime seems to be approximately quadratic in the
number of computed digits.

\section{Soundness}
\label{S:Sound}
We have extracted terms from the two proofs of Theorem~\ref{T:CoIDiv}
and Corollary~\ref{C:CoGDiv} and tested them on numerical examples.
Now we want to prove that these terms are correct, in the sense that
they \inquotes{satisfy their specification}.  We interpret this
statement in the sense of Kolmogorov \cite{Kolmogorov32}: a
(c.r.)\ formula should be viewed as a problem asking for a solution.
But what is a solution of a formula/problem $A$?  We use Kreisel's
notion of (modified) realizability and understand \inquotes{the term
  $t$ is a solution of $A$} as $t \mr A$ ($t$ realizes $A$).  This
definition is reviewed in Section~\ref{SS:Realizers}.  Then we go on
and prove in Section~\ref{SS:Sound} that for every proof $M$ of a
c.r.\ formula $A$ its extracted term $\extrTer{M}$ realizes $A$, i.e.,
$\extrTer{M} \mr A$.  In this proof we make use of
\inquotes{invariance axioms}\footnote{They are called (A{-}r)
  \inquotes{to assert is to realize} in \cite{Feferman79}.} stating
that every c.r.\ formula is invariant under realizability, formally $A
\leftrightarrow \ex_z(z \mr A)$.  Finally in Section~\ref{SS:FSound}
we report on a Minlog tool automatically generating such soundness
proofs.

\subsection{Realizers}
\label{SS:Realizers}
Recall that computational content arises from (co)inductive predicates
only.  Therefore we begin with a definition of what it means to be a
realizer of such a predicate.  We restrict ourselves to define
realizers for special instances only; a general treatment can be found
in \cite[p.334]{SchwichtenbergWainer12}.

Consider the definition of the inductive predicate $I_0$ in
Section~\ref{SS:Reals}.  By another (n.c.)\ inductive predicate
$I_0^\mrind$ of arity $(\typeR, \typeList)$ we can express that a list
$u$ witnesses (\inquotes{realizes}) that the real $x$ is in $I_0$.  We
write $u \mr I_0 x$ ($u$ is a realizer of $x \in I_0$) for $(x,u) \in
I_0^{\mrind}$.  The predicate $I_0^{\mrind}$ is required to be
non-computational, since in $(x,u) \in I_0^{\mrind}$ we already have a
realizer $u$.  $I_0^{\mrind}$ is inductively defined by the two
clauses
\begin{equation*}
  (0, \nil) \in I_0^{\mrind}, \quad
  \forall_{d \in \mathrm{Sd}} \forall_{(x,u) \in I_0^{\mrind}}\Bigl(
    \Bigl(\frac{x+d}{2}, s_d::u \Bigr) \in I_0^{\mrind} \Bigr)
\end{equation*}
and the induction axiom
\begin{equation*}
  (0, \nil) \in Q \to
  \forall_{d \in \mathrm{Sd}} \forall_{x \in I_0^{\mrind} \cap Q}\Bigl(
  \Bigl(\frac{x+d}{2}, s_d::u \Bigr) \in Q \Bigr) \to
  I_0^{\mrind} \subseteq Q.
\end{equation*}
where $s_d$ is the signed digit corresponding to $d \in \mathrm{Sd}$
\cite[p.334]{SchwichtenbergWainer12}).  We write $u \mr
I_0 x$ ($u$ is a realizer of $x \in I_0$) for $(x,u) \in
I_0^{\mrind}$.  The predicate $I_0^{\mrind}$ is required to be
non-computational, since in $(x,u) \in I_0^{\mrind}$ we already have a
realizer $u$.  $I_0^{\mrind}$ is inductively defined by the two
clauses
\begin{equation*}
  (0, \nil) \in I_0^{\mrind}, \quad
  \forall_{d \in \mathrm{Sd}} \forall_{(x,u) \in I_0^{\mrind}}\Bigl(
    \Bigl(\frac{x+d}{2}, s_d::u \Bigr) \in I_0^{\mrind} \Bigr)
\end{equation*}
and the induction axiom
\begin{equation*}
  (0, \nil) \in Q \to
  \forall_{d \in \mathrm{Sd}} \forall_{x \in I_0^{\mrind} \cap Q}\Bigl(
  \Bigl(\frac{x+d}{2}, s_d::u \Bigr) \in Q \Bigr) \to
  I_0^{\mrind} \subseteq Q.
\end{equation*}
where $s_d$ is the signed digit corresponding to $d \in \mathrm{Sd}$.
Similarly we coinductively define the n.c.\ predicate
$(\CoI_0)^{\mrind}$ of arity $(\typeR, \typeList)$ to express that a
list $u$ witnesses (\inquotes{realizes}) that the real $x$ is in
$\CoI_0$.  We write $u \mr \CoI_0 x$ ($u$ is a realizer of $x \in
\CoI_0$) for $(x,u) \in \CoI_0^{\mrind}$.  The closure axiom is
\begin{equation*}
  \forall_{(x,u) \in (\CoI_0)^{\mrind}}\Bigl( (x{=}0 \land u{=}\nil) \lor
 \ex_{d \in \mathrm{Sd}} 
  \ex_{(x',u') \in (\CoI_0)^{\mrind}}\Bigl(
   x=\frac{x'+d}{2} \land u=s_d::u' \Bigr) \Bigr)
\end{equation*}
and the coinduction axiom
\begin{align*}
  \forall_{(x,u) \in Q}\Bigl( &(x=0 \land u=\nil) \lor {}
  \\
  &\ex_{d \in \mathrm{Sd}} 
  \ex_{(x',u') \in (\CoI_0)^{\mrind} \cup Q}\Bigl(
   x=\frac{x'+d}{2} \land u=s_d::u' \Bigr) \Bigr) \to
   Q \subseteq (\CoI_0)^{\mrind}.
\end{align*}

\subsection{Soundness theorem}
\label{SS:Sound}
A formula or predicate is called $\mr$-free if it does not contain
any of the $I^\mrind$- or $\CoI^\mrind$-predicates.  A proof $M$ is
called $\mr$-free if it contains $\mr$-free formulas only.

\begin{thm}[Soundness]
  \label{T:Soundness}
  Let $M$ be an $\mr$-free proof of a formula $A$ from assumptions
  $u_i \colon C_i$ (\,$i<n$).  Then we can derive
  \begin{equation*}
    \begin{cases}
      \extrTer{M} \mr A &\hbox{if $A$ is c.r.}
      \\
      A &\hbox{if $A$ is n.c.}
    \end{cases}
  \end{equation*}
  from assumptions
  \begin{equation*}
    \begin{cases}
      z_{u_i} \mr C_i &\hbox{if $C_i$ is c.r.}
      \\
      C_i &\hbox{if $C_i$ is n.c.}
    \end{cases}
  \end{equation*}
\end{thm}

\begin{proof}
  By induction on $M$.  In the base case we have to prove that the
  extracted terms of $I^\pm$, $\CoI^\pm$ realize the respective
  axioms.  The proofs in \cite[Sections 7.2.8 and
    7.2.10]{SchwichtenbergWainer12} can be adapted to the present
  situation.  The step cases are easier; we only consider the ones
  where invariance axioms are used.

  \emph{Case} $(\lambda_{u^A} M^B)^{A \to B}$ with $A$ c.r.\ and $B$
  n.c.  We need a derivation of $A \to B$.  By induction hypothesis we
  have a derivation of $B$ from $z \mr A$.  Using the invariance axiom
  $A \to \ex_{z}( z \mr A)$ we obtain the required derivation of $B$
  from $A$ as follows.
  \begin{equation*}
    \AxiomC{$A \to \ex_z( z \mr A)$}
    \AxiomC{$A$}
    \BinaryInfC{$\ex_z( z \mr A)$}
    \AxiomC{$[z \mr A]$}
    \noLine
    \UnaryInfC{$\phantom{\hbox{IH}} \mid \hbox{IH}$}
    \noLine
    \UnaryInfC{$B$}
    \RightLabel{$\exE$}
    \BinaryInfC{$B$}
    \DisplayProof
  \end{equation*}

  \emph{Case} $(M^{A \to B}N^A)^B$ with $A$ c.r.\ and $B$ n.c.  The
  goal is to find a derivation of $B$.  By induction hypothesis we
  have derivations of $A \to B$ and of $\extrTer{N} \mr A$.  Now using
  the invariance axiom $\forall_z( z \mr A \to A)$ we obtain the
  required derivation of $B$ by $\to^-$ from the derivation of $A \to
  B$ and
  \begin{equation*}
    \AxiomC{$\forall_z( z \mr A \to A)$}
    \AxiomC{$\extrTer{N}$}
    \BinaryInfC{$\extrTer{N} \mr A \to A$}
    \AxiomC{$\phantom{\hbox{IH}} \mid \hbox{IH}$}
    \noLine
    \UnaryInfC{$\extrTer{N} \mr A$}
    \BinaryInfC{$A$}
    \DisplayProof
    \qedhere
  \end{equation*}
\end{proof}

\subsection{Formal soundness proofs for real division algorithms}
\label{SS:FSound}
In the files \texttt{sddiv.scm}, \texttt{graydiv.scm} the Minlog
command \texttt{add-sound} has been applied repeatedly to
c.r.\ theorems to automatically generate the corresponding soundness
proofs.

\section{Conclusion}
\label{S:Conclusion}
Recall that our goal was the extraction of computational content from
proofs in constructive analysis.  Although these proofs work with some
standard representation of real numbers, the technique used here makes
it possible to extract from such proofs terms describing algorithms
which operate on different representations of real numbers, for
instance streams of signed digits or Gray code.  Moreover, formal
proofs of their correctness can be generated automatically.

As future work it is planned to extend the present approach to
problems involving e.g.\ the exponential function and more generally
power series.  A promising application area would be Euler's existence
proof of solutions for ordinary differential equations satisfying a
Lipschitz condition.

\section*{Acknowledgment}
We thank Tatsuji Kawai, Nils K\"opp, Kenji Miyamoto, Hideki Tsuiki
and three anonymous referees for helpful comments.

\bibliographystyle{alpha}
\bibliography{div4}

\end{document}